\documentclass[12pt]{article}
\pdfoutput=1
\usepackage{graphicx}
\usepackage{listings}
\usepackage{mathrsfs}
\usepackage[font={small},labelfont={sf,bf}]{caption}
\usepackage{color}
\usepackage{amssymb}
\usepackage[utf8]{inputenc}
\usepackage{afterpage}
\usepackage[T1]{fontenc}
\usepackage{amsfonts}
\usepackage{mathtools}
\usepackage{amsmath}
\usepackage{amssymb}
\usepackage{verbatim}
\usepackage[margin=1.2in]{geometry}
\usepackage{bm}
\usepackage{bbm}
\usepackage{hyperref}
\usepackage{enumitem}
\usepackage{stmaryrd}
\usepackage{authblk}
\usepackage[numbers,sort]{natbib}

\usepackage{algorithm,algcompatible}
\usepackage{multirow}
\usepackage[numbers,sort]{natbib}

\graphicspath{{img/}}

\newcommand*\mcupinn[2]{\vcenter{\hbox{$\mathsurround=0pt
  \ifx\displaystyle#1\textstyle\else#1\fi\bigcup$}}}

\newcommand*\mcapinn[2]{\vcenter{\hbox{$\mathsurround=0pt
  \ifx\displaystyle#1\textstyle\else#1\fi\bigcap$}}}
\DeclarePairedDelimiter{\abs}{\lvert}{\rvert}
\DeclarePairedDelimiter{\norm}{\lVert}{\rVert}
\DeclarePairedDelimiter{\parr}{(}{)}
\DeclarePairedDelimiter{\parq}{[}{]}

\DeclarePairedDelimiter{\bra}{\lbrace}{\rbrace}

\DeclarePairedDelimiter{\prodscal}{\langle}{\rangle}

\DeclareMathOperator*{\argsup}{arg\,sup}
\DeclareMathOperator*{\expval}{\mathbb{E}}

\DeclareMathOperator*{\st}{\,:\,}
\usepackage{authblk}

\title{A Weighted POD Method for Elliptic PDEs with Random Inputs \thanks{This work was funded by European Union Funding for Research and Innovation (project H2020 ERC
CoG 2015 AROMA-CFD project 681447) and by the INDAM-GNCS project.}}
\author[1]{Luca Venturi \thanks{Current address: Courant Institute of Mathematical Sciences, New York University, New York, U.S. Corresponding author: venturi@cims.nyu.edu}}
\author[1]{Francesco Ballarin }
\author[1]{Gianluigi Rozza }
\affil[1]{mathLab, Mathematics Area, SISSA, Trieste, Italy}

\begin{document}
\maketitle

\begin{abstract}
In this work we propose and analyze a weighted proper orthogonal decomposition method to solve elliptic partial differential equations depending on random input data, for stochastic problems that can be transformed into parametric systems. 
The algorithm is introduced alongside the weighted greedy method. Our proposed method aims to minimize the error in a $L^2$ norm and, in contrast to the weighted greedy approach, it does not require the availability of an error bound.
Moreover, we consider sparse discretization of the input space in the construction of the reduced model; for high-dimensional problems, provided the sampling is done accordingly to the parameters distribution, this enables a sensible reduction of computational costs, while keeping a very good accuracy with respect to high fidelity solutions.
We provide many numerical tests to asses the performance of the proposed method compared to an equivalent reduced order model without weighting, as well as to the weighted greedy approach, in both low and higher dimensional problems.
\end{abstract}

\section{Introduction}
Many physical models and engineering problems are described by partial differential equations (PDEs). These equations usually depends on a number of parameters, which can correspond to physical or geometrical properties of the model as well as to initial or boundary conditions of the problem, and can be deterministic or subject to some source of uncertainty. Many applications require either to compute the solutions in real time, i.e. to compute them almost instantaneously given certain input parameters, or to compute a huge number of solutions corresponding to diverse parameters. The latter is for example the case when we want to evaluate statistics of some features associated to the solution. In both cases it is necessary to have methods to get very fast, and still reliable, numerical solutions to these equations. \\
Reduced order models \cite{quarteroni2011certified,rozza:book} were developed to face this problem. The main idea of these methods is that in many practical cases the solutions obtained for different values of the parameters belong to a lower-dimensional manifold of the functional space they naturally belong to. Therefore, instead of looking for the solution in a high-dimensional approximation subspace of this functional space, one could look for it in a low-dimensional, suitably tailored, \emph{reduced order space}. Obviously, this considerably reduces the computational effort, once the reduced spaces are available. Different methods have been developed to construct these sub-spaces in an efficient and reliable way: two popular techniques are the greedy method and the Proper Orthogonal Decomposition (POD) method. \\  
As said above, many of the parameters describing the problem at hand can be subject to uncertainties. In several applications, keeping track of these uncertainties could make the model more flexible and more reliable than a deterministic model. In this case we talk about PDEs with random inputs, or \emph{stochastic} PDEs.   
Many different methods are available to solve such problems numerically, such as stochastic Galerkin \cite{sullivan:uq}, stochastic collocation \cite{nobile:collocation,hestaven:collocation,lorenzo}, as well as Monte-Carlo \cite{sullivan:uq}. More recently, methods based on reduced order models have been developed \cite{peng:thesis,peng:w_elliptic_rnd,peng:comparison,chen2017reduced,spannring2017weighted,torlo2017stabilized}. In practice, they aim to construct low-dimensional approximation subspaces in a \emph{weighted} fashion, i.e. more weight is given to more relevant values of the parameters, according to an underlying probability distribution. Nevertheless, to our knowledge, no reduced order method for stochastic problems based completely on the POD algorithm has been proposed so far.
\\
In this paper we consider the case of linear elliptic PDEs with the same assumption as in \cite{peng:w_elliptic_rnd}. We briefly discuss the framework and introduce the weighted greedy algorithm. The rest of the work is devoted to define a weighted proper orthogonal decomposition method, which is the original contribution of this work. While the greedy algorithm aims to minimize the error in a $L^\infty$ norm over the parameter space, the POD method looks instead at the $L^2$ norm error, which, in the stochastic case, is the mean square error. The integral of the square error is discretized according to some quadrature rule, and the reduced order space is found by minimizing this approximated error. From the algorithmic point of view, this corresponds to applying a principal component analysis to a weighted (pre-conditioned) matrix. The main difference from the greedy algorithm is that our method does not require the availability of an upper bound on the error. Nevertheless, our method requires an \emph{Offline} computation of an higher number of high-fidelity solutions, which may be costly. Therefore, to keep the computational cost low, we adopted sparse quadrature rules. These methods allow to get the same accuracy at a sensibly lower computational cost, especially in the case of high dimensional parameter space. The resulting weighted POD method shows the same accuracy as long as the underlying quadrature rule is properly chosen, i.e. accordingly to the distribution. Indeed, since the reduced order space will be a subspace of the space spanned by the solutions on the parameter nodes given by the quadrature rule, these training parameters needs to be representative of the probability distribution on the parameter space. Thus the choice of a proper quadrature rule plays a fundamental role. 
\\
The work is organized as it follows. An elliptic PDE with random input data is
set up with appropriate assumptions on both the random coefficient and the forcing
term in Section 2. 
In Section 3 we describe the affinity assumptions that must be made on our model to allow the use of the reduced basis methods; the greedy algorithm is presented in its weighted version. Our proposed method is then described as based on the minimization of an approximated mean square error. Sparse interpolation techniques for a more efficient integral approximation in higher dimensions are also presented. 
Numerical examples (Section 4) for both a low-dimensional ($N=4$) problem and its higher dimensional counterpart ($N=9$) are presented as verification of the efficiency and convergence properties of our method.
Some brief concluding remarks are made in Section 5.

\section{Problem setting}

Let $(\Omega,\mathcal{F},P)$ be a complete probability space, where $\Omega$ is the set of possible outcomes, $\mathcal{F}$ is the $\sigma$-algebra of events and $P$ is the probability measure. Moreover, let $D\subseteq\mathbb{R}^d$ ($d=1,2,3$) be an open, bounded and convex domain with Lipschitz boundary $\partial D$. We consider the following stochastic elliptic problem: find $u:\Omega \times \overline{D} \to \mathbb{R}$, such that it holds for $P$-a.e. $\omega\in\Omega$ that
\begin{align}
    - \nabla \cdot (a(\omega,x)\cdot\nabla u(\omega, x)) & = f(\omega,x), \qquad x \in D, \label{eq:stoc_ell_pde} \\
    u(\omega,x) & = 0, \qquad \in \partial D. \notag
\end{align}
Here $a$ is a strictly positive random diffusion coefficient and $f$ is a random forcing term; the operator $\nabla$ is considered w.r.t. the spatial variable $x$. 
Motivated by the regularity results obtained in \cite{peng:w_elliptic_rnd}, we make the following assumptions:
\begin{enumerate}
    \item The forcing term $f$ is square integrable in both variables, i.e.
    $$
    \norm{f}^2_{L^2_P(\Omega)\otimes L^2(D)} \doteq \int_{\Omega \times D} f^2(\omega,x)\,dP(\omega)dx < \infty.
    $$
    \item The diffusion term is $P$-a.s. uniformly bounded, i.e. there exist $0<a_{\text{min}}<a_{\text{max}}<\infty$ such that
    $$
    P\parr*{a(\cdot,x)\in (a_{\text{min}},a_{\text{max}}) \,\forall\, x \in D}=1.
    $$
\end{enumerate}
If we introduce the Hilbert space $\mathbb{H} = L^2_P(\Omega)\otimes H^1_0(D)$, we can consider the weak formulation of problem \eqref{eq:stoc_ell_pde}: find $u\in \mathbb{H}$ s.t.
\begin{equation}\label{eq:stoc_ell_pde_ww}
\mathbb{A}(u,v) = \mathbb{F}(v), \qquad \text{for every }v\in\mathbb{H},
\end{equation}
where $\mathbb{A}:\mathbb{H}\times\mathbb{H}\to\mathbb{R}$ and $\mathbb{F}:\mathbb{H}\to\mathbb{R}$ are,  respectively, the bilinear and linear forms
$$
\mathbb{A}(u,v) = \int_{\Omega\times D} a\, \nabla u \cdot \nabla v \,dP(\omega) dx, \qquad \mathbb{F}(v) = \int_{\Omega\times D} v\cdot f\,dP(\omega) dx.
$$
The above is called weak-weak formulation. Thanks to assumption i.-ii., the Lax-Milgram theorem \cite{brezis} ensures us the existence of a unique solution $u\in\mathbb{H}$.

More than the solution itself, we will be interested in statistics of values related to the solution, e.g., the expectation
$\expval[s(u)]$, where $s(u)$ is some linear functional of the solution. In particular, the numerical experiments in the following are all performed for the computation of $\expval\parq{u}$. 

\subsection{Weak-strong formulation}
  
To numerically solve problem \eqref{eq:stoc_ell_pde_ww} we first need to reduce $(\Omega,\mathcal{F},P)$ to a finite dimensional probability space. This can be accomplished up to a desired accuracy through, for example, the Karhunen-Lo\'eve expansion \cite{loeve:probability}. The random input is in this case characterized by $K$ uncorrelated random variables $Y_1(\omega),\dots,Y_K(\omega)$ so that we can write
$$
a(\omega,\cdot) = a(\mathbf{Y}(\omega),\cdot), \qquad
f(\omega,\cdot) = f(\mathbf{Y}(\omega),\cdot),
$$
and hence (thanks to the Doob-Dynkin lemma \cite{oksendal})
$$
u(\omega,\cdot) = u(\mathbf{Y}(\omega),\cdot),
$$
where $\mathbf{Y}(\omega) = (Y_1(\omega),\dots,Y_K(\omega))$. We furthermore assume that $\mathbf{Y}$ has a continuous probability distribution with density function $\rho:\mathbb{R}^K\to\mathbb{R}^+$ and that $Y_k(\Omega) \subseteq \Gamma_k$ for some $\Gamma_k\subset \mathbb{R}$ compact sets, $k=1,\dots,K$.  
In case the initial probability density is not compactly supported, we can easily reduce to this case by truncating the function $\rho$ on a compact set up to a desired accuracy. Our problem can be reformulated in terms of a weighted parameter $y\in\Gamma\doteq \prod_{k=1}^K\Gamma_k$: find $u:\Gamma\to \mathbb{V} \doteq H_0^1(D)$ such that
\begin{equation}\label{eq:stoc_ell_pde_ws}
A(u(y),v;y) = F(v;y) \qquad \text{for all } v\in\mathbb{V},
\end{equation}
for a.e. $y\in \Gamma$, where $A(\cdot,\cdot;y)$ and $F(\cdot;y)$ are, respectively, the parametrized bilinear and linear forms defined as
$$
A(u,v;y) = \int_D a(y,x)\, \nabla u(x) \cdot \nabla v(x) \,dx, \qquad F(v;y) = \int_D v(x) f(y,x)\,dx,
$$
for $k=1,\dots,K$.
The parameter $y$ is distributed according to the probability measure $\rho(y)dy$. Problem \eqref{eq:stoc_ell_pde_ws} is called the weak-strong formulation of problem \eqref{eq:stoc_ell_pde_ww}. Again, the existence of a solution is guaranteed by the Lax-Milgram theorem. 

Given an approximation space $\mathbb{V}_\delta \subseteq \mathbb{V}$ (typically a finite element space), with $\mathrm{dim}(\mathbb{V}_\delta) = N_\delta < \infty$, we consider the approximate problem: find $u_{N_\delta}:\Gamma \to \mathbb{V}_\delta$ such that
\begin{equation}\label{eq:stoc_ell_pde_fe}
A(u_{N_\delta}(y),v;y) = F(v;y) \qquad \text{for all } v\in\mathbb{V_\delta},
\end{equation}
for a.e. $y\in \Gamma$. We refer to problem \eqref{eq:stoc_ell_pde_fe} as the \emph{truth} (or high dimensional) problem and $u_{N_\delta}$ as the \emph{truth} solution. Consequently we approximate the output of interest as $s_{N_\delta} = s(u_{N_\delta}) \simeq s(u)$ and its statistics, e.g. $\expval[s_{N_\delta}] \simeq \expval[s(u)]$.
  
\subsection{Monte-Carlo methods}

A typical way to numerically solve problem \eqref{eq:stoc_ell_pde_ws} is to use a Monte-Carlo simulation. This procedure takes the following steps:
\begin{enumerate}
\item generate $M$ (given number of samples) independent and identically distributed (i.i.d.) copies of $\mathbf{Y}$, and store the obtained values $y^j$, $j=1,\dots,M$; 
\item solve the deterministic problem \eqref{eq:stoc_ell_pde_ws} with $y^j$ and obtain solution $u_j=u(y^j)$, for $j=1,\dots,M$;
\item evaluate the solution statistics as averages, e.g.,
$$
\mathbb{E}[u] \simeq \langle u \rangle  = \frac{1}{M}\sum_{j=1}^M u_j \qquad \textrm{or} \qquad
\mathbb{E}[s(u)] \simeq \langle s(u) \rangle  = \frac{1}{M}\sum_{j=1}^M s(u_j) 
$$
for some suitable output function $s$.
\end{enumerate}
Although the convergence rate of a Monte Carlo method is formally independent of the dimension of the random space, it is relatively slow (typically $1/\sqrt{M}$). Thus, one requires to solve a large amount of deterministic problems to obtain a desired accuracy, which implies a very high computational cost. In this framework reduced order methods turn out to be very useful in order to reduce the computational cost, at cost of a (possibly) small additional error.

\section{Weighted reduced order methods}\label{section:w_rom}

Given the truth approximation space $\mathbb{V}_\delta$, reduced order algorithms look for an approximate solution of \eqref{eq:stoc_ell_pde_fe} in a reduced order space $\mathbb{V}_N\subseteq\mathbb{V}_\delta$, with $N\ll N_\delta$. Reduced order methods consist of an \emph{offline} phase and an \emph{online} phase. 

During the offline phase, an hierarchical sequence of reduced order spaces $\mathbb{V}_1\subseteq\cdots\subseteq\mathbb{V}_{N_{\mathrm{max}}}$ is built. These spaces are sought in the subspace spanned by the solutions for a discrete training set of parameters $\Xi_t = \bra{y^1,\dots,y^{n_t}}\subseteq \Gamma$,
according to some specified  algorithm.   
In this work we focus on the greedy and POD algorithm, as we will detail in Sections 3.2 and 3.3, respectively. 

During the online phase, the reduced order problem is solved: find $u_N:\Gamma \to \mathbb{V}_N$ such that
\begin{equation}\label{eq:stoc_ell_pde_ro}
A(u_N(y),v;y) = F(v;y) \qquad \text{for all } v\in\mathbb{V}_N,
\end{equation}
for a.e. $y\in \Gamma$. At this point, we can approximate the output of interest as $s_N = s(u_N) \simeq s(u)$ and its statistics, e.g. $\expval[s_N] \simeq \expval[s(u)]$. In the stochastic case, we would also like to take into account the probability distribution of the parameter $y\in\Gamma$. \emph{Weighted} reduced order methods consist of  slight modifications of the offline phase so that a \emph{weight} $w(y_i)$ is associated to each sample parameter $y_i$, according to the probability distribution $\rho(y)dy$. As we will discuss, another crucial point is the choice of the training set $\Xi_t$.

For sake of simplicity, from now on we omit the indexes $\delta$, $N_\delta$ and we assume that our original problem \eqref{eq:stoc_ell_pde_ws} coincides with the truth problem \eqref{eq:stoc_ell_pde_fe}.

\subsection{Affine Decomposition assumption}

In order to ensure efficiency of the online evaluation, we need to assume that the bilinear form $A(\cdot,\cdot;y)$ and linear form $F(\cdot;y)$ admit an \emph{affine decomposition}. In particular, we require the diffusion term $a(y,x)$ and the forcing term $f(y,x)$ to have the following form:
\begin{align*}
    a(x,y) & = a_0(x) + \sum_{k=1}^K a_k(x)y_k, \\
    f(x,y) & = f_0(x) + \sum_{k=1}^K f_k(x)y_k,
\end{align*}
where $a_k\in L^\infty(D)$ and $f_k\in L^2(D)$, for $k=1,\dots,K$, and $y = (y_1,\dots,y_K)\in\Gamma$. Thus, the bilinear form $A(\cdot,\cdot;y)$ and linear form $F(\cdot;y)$ can be written as 
\begin{equation}
A(u,v;y) = A_0(u,v) + \sum_{k=1}^K y_k A_k(u,v), \qquad F(u,v;y) = F_0(u,v) + \sum_{k=1}^K y_k F_k(u,v) \label{eq:aff_dec}
\end{equation}
with
$$
A_k(u,v) = \int_D a_k\, \nabla u \cdot \nabla v \,dx, \qquad F_k(v) = \int_D v f_k\,dx,
$$
for $k=1,\dots,K$. In the more general case that the functions $a(\cdot; y)$, $f(\cdot; y)$ do not depend on $y$ linearly, one can reduce to this case by employing the empirical interpolation method \cite{patera:eim}. A weighted version of this algorithm has also been proposed in \cite{peng:eim}.

\subsection{Weighted greedy algorithm}

The greedy algorithms ideally aims to find the $N$-dimensional subspace that minimizes the error in the $L^\infty(\Gamma)$ norm, i.e. $\mathbb{V}_N\subseteq \mathbb{V}$ such that the quantity
\begin{equation}
\sup_{y\in\Gamma} \norm{u(y)-u_N(y)}_\mathbb{V}. \label{eq:gre_min_qua}
\end{equation}
is minimized.
The reduced order spaces are build hierarchically as 
$$
\mathbb{V}_N = \mathrm{span}\bra{u(y^1),\dots,u(y^N)},
$$
for $N=1,\dots,N_{\mathrm{max}}$, where the parameters $y^N$ are sought as solution of an $L^\infty$ optimization problem in a greedy way:
$$
y^N = \argsup_{y\in\Gamma} \norm{u(y)-u_{N-1}(y)}_\mathbb{V}.
$$
To actually make the optimization problem above computable, the parameter domain $\Gamma$ is replaced with the training set $\Xi_t$. The greedy scheme is strongly based on the availability of an error estimator $\eta_N(y)$ such that
$$
\norm{u(y)-u_N(y)}_\mathbb{V} \leq \eta_N(y).
$$
Such an estimator needs to be both \emph{sharp}, meaning that there exists a constant $c$ (possibly depending on $N$) as close as possible to $1$ such that $\norm{u(y)-u_N(y)}_\mathbb{V} \geq c\cdot\eta_N(y)$, and \emph{efficiently} computable. An efficient way to compute $\eta_N$ is described in detail e.g. in \cite{rozza:book} (the affine decomposition \eqref{eq:aff_dec} plays an essential role here). The greedy scheme actually looks for the maximum of the estimator $\eta_N$ instead of the quantity in \eqref{eq:gre_min_qua} itself. In the case of stochastic parameter $y\in\Gamma$, the idea is to modify $\eta_N(y)$, multiplying it by a weight $w(y)$, chosen accordingly to the distribution $\rho(y)dy$. The estimator $\eta_N(y)$ is thus replaced by $\widehat{\eta}_N(y) \doteq w(y)\eta_N(y)$ in the so called weighted scheme. Weighted greedy methods have been originally proposed and developed in \cite{peng:thesis,peng:w_elliptic_rnd,peng:comparison,chen2017reduced,peng:multilevel}. An outline of the weighted greedy algorithm for the construction of the reduced order spaces is reported in Algorithm \ref{alg:w_gre}.
\begin{algorithm}[tb]
\caption{Weighted greedy algorithm}\label{alg:w_gre}
\begin{algorithmic}[1]
    \STATE sample a training set $\Xi_t \subseteq \Gamma$;
    \STATE choose the first sample parameter $y^1\in\Xi_t$;
    \STATE solve problem \eqref{eq:stoc_ell_pde_fe} at $y^1$ and construct $\mathbb{V}_1 = \mathrm{span}\bra{u(y^1)}$;
    \FOR{$N=2,\dots,N_{\mathrm{max}}$} 
        \STATE{compute the weighted error estimator $\widehat{\eta}_N(y)$, for all $y\in\Xi_t$;} 
        \STATE{choose $y^N = \argsup_{y\in\Xi_t} \widehat{\eta}_N(y)$;}
        \STATE{solve problem \eqref{eq:stoc_ell_pde_fe} at $y^N$ and construct $\mathbb{V}_N = \mathbb{V}_{N-1}\oplus\mathrm{span}\bra{u(y^N)}$;}
    \ENDFOR
\end{algorithmic}
\end{algorithm}
A greedy routine with estimator $\widehat{\eta}_N$ hence aims to minimize the distance of the solution manifold from the reduced order space in a weighted $L^\infty$ norm on the parameter space. Now, depending on what one wants to compute, different choices can be made for the weight function $w$. For example, if we are interested in a statistics of the solution, e.g., $\mathbb{E}[u]$, we can choose $w(y)=\rho(y)$. Thus, for the error committed computing the expected value using the reduced basis, the following estimates holds:
$$
\big\lVert \mathbb{E}[u]-\mathbb{E}[u_N] \big\rVert_\mathbb{V} \leq \int_\Gamma \norm{u(y)-u_N(y)}_\mathbb{V} \rho(y) \, dy
\leq \abs{\Gamma} \sup_{y\in\Gamma}\widehat{\eta}_N(y).
$$
If we are interested in evaluating the expectation of a linear output $s(u)$, $\mathbb{E}[s(u)]$, using the same weight function, we get the error estimate:
$$
\big| \mathbb{E}[s(u)]-\mathbb{E}[s(u_N)]\big|  \leq \int_\Gamma \norm{s}_\mathbb{V'}\norm{u(y)-u_N(y)}_\mathbb{V} \rho(y) \, dy
\leq \norm{s}_\mathbb{V'}\cdot\abs{\Gamma}\sup_{y\in\Gamma}\widehat{\eta}_N(y).
$$
Instead, taking $w(y)=\sqrt{\rho(y)}$ one gets the estimate for the quadratic error:
\begin{equation}
\label{eq:w_gre_sqrt_est}
\norm{u(\mathbf{Y})-u_N(\mathbf{Y})}_\mathbb{H}^2 = \expval\norm{u-u_N}_\mathbb{V}^2 \leq \abs{\Gamma} \sup_{y\in\Gamma} \widehat{\eta}_N(y)^2.
\end{equation} 

\subsection{Weighted POD algorithm}

The Proper Orthogonal Decomposition is a different method to build reduced order spaces, which does not require the evaluation of an error bound. The main idea in the deterministic (or uniform, i.e. when $\rho\equiv 1$) case is to find the $N$-dimensional subspace that minimizes the error in the $L^2(\Gamma)$ norm, i.e. $\mathbb{V}_N\subseteq \mathbb{V}$ such that the quantity
\begin{equation}
\int_\Gamma \norm{u(y)-u_N(y)}_\mathbb{V}^2\,dy. \label{eq:pod_min_qua}
\end{equation}
is minimized. As with the greedy algorithm, the method does not aim to minimize the quantity \eqref{eq:pod_min_qua} directly, but a discretization of it, namely
$$
\sum_{y\in\Xi_t}\norm{u(y)-u_N(y)}_\mathbb{V}^2 = \sum_{i=1}^{n_t}\norm{\varphi_i-P_N(\varphi_i)}_\mathbb{V}^2, 
$$
where $\varphi_i = u(y^i)$, $i=1,\dots,n_t$, and $P_N:\mathbb{V}\to\mathbb{V}_N$ the projection operator associated with the subspace $\mathbb{V}_N$. One can show (see e.g. \cite{rozza:book}) that the minimizing $N$-dimensional subspace $\mathbb{V}_N$ is given by the subspace spanned by the $N$ leading eigenvectors of the linear map
$$
c : \mathbb{V} \to \mathbb{V}, \quad v \mapsto c(v) = \sum_{i=1}^{n_t} \prodscal{v,\varphi_i}_\mathbb{V}\cdot \varphi_i,
$$
where $\prodscal{\cdot, \cdot}_\mathbb{V}$ denotes the inner product in $\mathbb{V}$.
Computationally, this is equivalent to find the $N$ leading eigenvectors of the symmetric matrix $C\in\mathbb{R}^{n_t\times n_t}$ defined as $C_{ij} = \prodscal{\varphi_i,\varphi_j}_\mathbb{V}$. In the case of stochastic inputs we would rather like to find the $N$-dimensional subspace that minimizes the following error 
\begin{equation}
\norm{u(\mathbf{Y})-u_N(\mathbf{Y})}_\mathbb{H}^2 = \int_\Gamma \norm{u(y)-u_N(y)}_\mathbb{V}^2\rho(y)\,dy. 
\label{eq:w_pod_min_qua}
\end{equation}
Based on this observation, we propose a weighted POD method, which is based on the minimization of a discretized version of \eqref{eq:w_pod_min_qua}, namely
\begin{equation}
\sum_{y\in\Xi_t}w(y)\norm{u(y)-u_N(y)}_\mathbb{V}^2 = \sum_{i=1}^{n_t}w_i\norm{\varphi_i-P_N(\varphi_i)}_\mathbb{V}^2, 
\label{eq:w_pod_min_qua_dis}
\end{equation}
where $w:\Xi_t \to [0,\infty)$ is a weight function prescribed according to the parameter distribution $\rho(y)dy$, and $w_i=w(y^i)$, $i=1,\dots,n_t$. Again, this is computationally equivalent to finding the $N$ maximum eigenvectors of the preconditioned matrix $\hat{C} \doteq P \cdot C$, where $C$ is the same as defined before and $P=\mathrm{diag}(w_1,\dots,w_{n_t})$. We note that the matrix $\hat{C}$ is not symmetric in the usual sense, but it is with respect to the scalar product induced by the matrix $C$ (i.e. it holds that $\hat{C}^TC = C\hat{C}$).
Thus, spectral theorem still holds and there exists an orthonormal basis of eigenvectors, i.e., $\widehat{C}$ is diagonalizable with an orthogonal change of basis matrix. The discretized parameter space $\Xi_t$ can be selected with a sampling technique, e.g., using an equispaced tensor product grid on $\Gamma$ or taking $M$ realizations of a uniform distribution on $\Gamma$. Note that if we build $\Xi_t$ as the set of $M$ realizations of a random variable on $\Gamma$ with distribution $\rho(y)dy$, and we put $w\equiv 1$, the quantity \eqref{eq:w_pod_min_qua_dis} we minimize is just a Monte Carlo approximation of the integral \eqref{eq:w_pod_min_qua}. Following this observation, a possible approach would be to select $\Xi_t$ and $w$ as the nodes and the weights of a quadrature rule that approximates the integral \eqref{eq:w_pod_min_qua}.
That is, if we consider a quadrature operator $\mathcal{U}$, defined as 
$$
\mathcal{U}(f) = \sum_{i=1}^{n_t} \omega_i f(x^i)
$$
for every integrable function $f:\Gamma\to\mathbb{R}$, where $\bra{x^1,\dots,x^{n_t}}\subset\Gamma$ are the nodes of the algorithm and $\omega_1,\dots,\omega_{n_t}$ the respective weights, then we can take $\Xi_t = \bra{x^1,\dots,x^{n_t}}$ and%
\footnote{We assume that $\mathcal{U}$ is a quadrature rule for integration with respect to $dy$. If a quadrature rule $\mathcal{U}_{\rho}$ for integration with respect to the weighted measure $\rho dy$ is used instead, that is suffices to take $w_i = \omega_i$.}
$w_i = \omega_i\rho(x^i)$.
Therefore, varying the quadrature rule $\mathcal{U}$ used, one obtain diverse ways of (at the same time) sampling the parameter space and preconditioning the matrix $C$. An outline of the weighted POD algorithm for the construction of the reduced order spaces is reported in Algorithm \ref{alg:w_pod}. \\
\begin{algorithm}[tb]
\caption{Weighted POD algorithm}\label{alg:w_pod}
\begin{algorithmic}[1]
    \STATE choose a training set $\Xi_t = \bra{y^1,\dots,y^{n_t}} \subseteq \Gamma$ and a weight function $w$ according to some quadrature method;
    \STATE compute the solutions $\varphi_i$ by solving problem \eqref{eq:stoc_ell_pde_fe} at $y^i$, $i=1,\dots,n_t$; 
    \STATE assemble the matrix $\hat{C}_{ij} = w_i\prodscal{\varphi_i,\varphi_j}_\mathbb{V}$ and compute its $N$ maximum eigenvectors $n^1,\dots,n^j$;
    \STATE compute $N$ maximal eigenvectors $\xi^1,\dots,\xi^N\in\mathbb{V}$ as $\xi^i = \sum_{j=1}^{n_t} n^i_j \varphi_j$ and construct $\mathbb{V}_N = \mathrm{span}\bra{\xi^1,\dots,\xi^N}$;
\end{algorithmic}
\end{algorithm}
Since the proposed weighted POD method requires to perform $n_t$ truth solve and compute the eigenvalues of a $n_t\times n_t$ matrix, the dimension $n_t$ of $\Xi_t$ should better not be too large. A possible way to keep the sample size $n_t$ low is to adopt a sparse grid quadrature rule, as we describe in the following.

\subsection{Sparse grid interpolation}

In order to make the weighted POD method more efficient as the dimension of the parameter space increases, we use Smolyak type sparse grid instead of full tensor product ones for the quadrature operator $\mathcal{U}$. These type of interpolation grids have already been used in the context of weighted reduced order methods \cite{peng:thesis,peng:comparison} as well as of other numerical methods for stochastic partial differential equations, like e.g. in stochastic collocation \cite{nobile:collocation,hestaven:collocation}. Full tensor product quadrature operator are simply constructed as product of univariate quadrature rules. If $\Gamma = \prod_{k=1}^K \Gamma_k$ and $\bra{\mathcal{U}_i^{(k)}}_{i=1}^\infty$ are sequences of univariate quadrature rules on $\Gamma_k$, $k=1,\dots,K$, the tensor product multivariate rule on $\Gamma$ of order $q$ is given by
$$
\mathcal{U}^K_q \doteq \bigotimes_{k=1}^K\,\mathcal{U}_q^{(k)} = \sum_{\substack{\abs{\alpha}_\infty \leq q \\ \alpha \in \mathbb{N}^K}} \bigotimes_{k=1}^K \Delta_{\alpha_k}^{(k)},
$$
where we introduced the \emph{differences operators} $\Delta_0^{(k)}=0$, $\Delta_{i+1}^{(k)}=\mathcal{U}_{i+1}^{(k)}-\mathcal{U}_i^{(k)}$ for $i\geq 0$. Given this, the Smolyak quadrature rule of order $q$ is defined as
$$
\mathcal{Q}_q^K \doteq \sum_{\substack{\abs{\alpha}_1 \leq q \\ \alpha \in \mathbb{N}^K}} \bigotimes_{k=1}^K \Delta_{\alpha_k}^{(k)}.
$$
Therefore, Smolyak type rules can be seen as delayed sum of ordinary full tensor product rules. One of their main advantage is that the number of evaluation points is drastically reduced as $K$ increases. Indeed, if $X^{(k)}_i$ is the set of evaluation points for the rule $\mathcal{U}_i^{(k)}$, then the set of evaluation points for the rule $\mathcal{U}_q^K$ is given by
$$
\Theta_F^{q,K} =  X_q^{(1)} \times \dots \times X_q^{(K)},
$$
while the set of evaluation points for the rule $\mathcal{Q}_q^K$ is given by
$$
\Theta_S^{q,K} = \bigcup_{\substack{\abs{\alpha}_1 = q \\ \alpha \in \mathbb{N}^n, \; \alpha \geq \mathbf{1}}} X_{\alpha_1}^{(1)} \times \dots \times X_{\alpha_n}^{(n)}.
$$
In practice, for $K$ large, one has that $\abs{\Theta_S^{q,K}}$ grows much slower as $q$ increases w.r.t. $\abs{\Theta_F^{q,K}}$, which has an exponential behavior (see Fig. 1 for a comparison using a Clenshaw-Curtis quadrature rule). This implies that the adoption of Smolyak quadrature rules has a much lower computational cost than common full tensor product rules for high dimensional problem. Moreover, the performance of Smolyak quadrature rules are comparable with standard rules. For more details about error estimates and computational cost of Smolyak rules we refer e.g. to \cite{novak:sparse_grids,gerstner:sparse_grids,holtz:sparse_grids,novak:sparse_grids_quadrature,wasil:sparse_grids}.
\begin{figure}[tb]
\centering
\begin{minipage}[c]{.47\textwidth}
\includegraphics[width=\textwidth,
keepaspectratio]{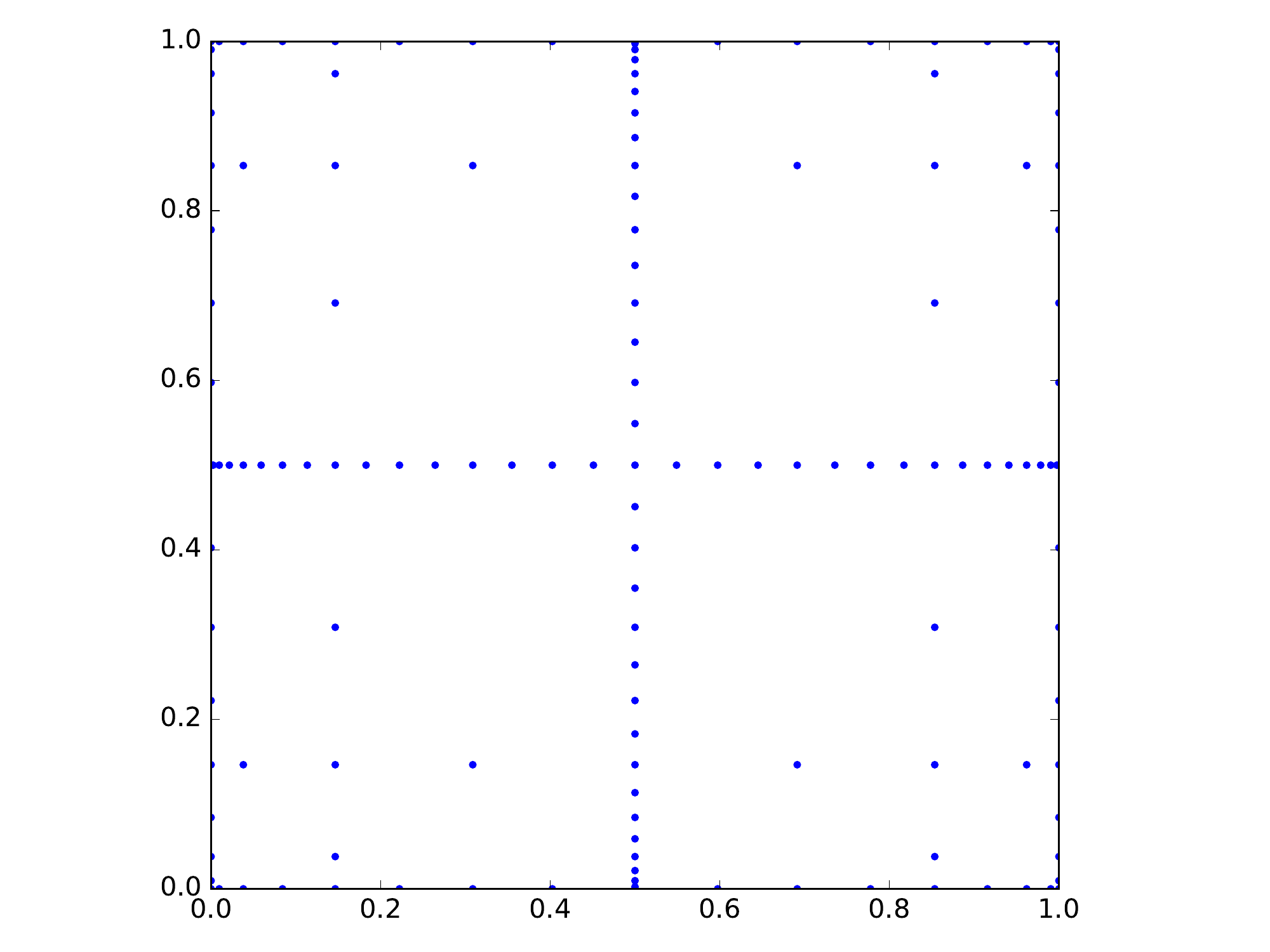}
\end{minipage}
\hspace{4mm}
\begin{minipage}[c]{.47\textwidth}
\includegraphics[width=\textwidth,
keepaspectratio]{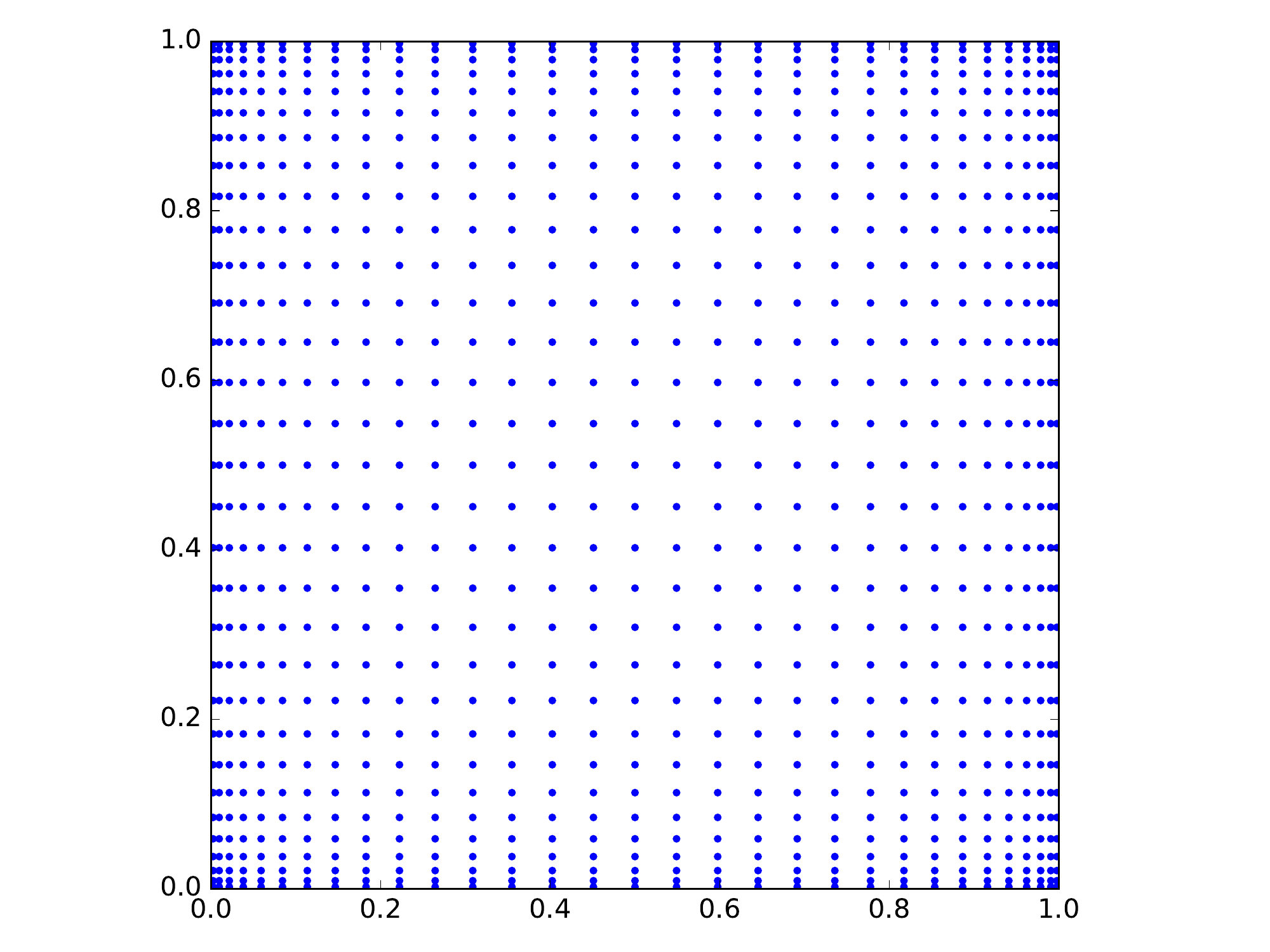}
\end{minipage}
\caption{Two dimensional grids based on nested Clenshaw-Curtis nodes for $q=6$. The left one is based on a Smolyak rule ($145$ nodes), while the right one on a full tensor product rule ($1089$ nodes).}
\end{figure}

\section{Numerical tests}

We tested and compared the weighted greedy algorithm and the weighted POD algorithm for the solution of problem \eqref{eq:stoc_ell_pde} on the unit square domain $D = [0,1]^2$. To solve the below problems we used the $\mathtt{RBniCS}$ library \cite{rbnics:link}, built on top of $\mathtt{FEniCS}$ \cite{logg2012automated}. 

\subsection{Thermal block problem}\label{section:4dz}

\begin{figure}[tb]
\centering
\begin{minipage}[c]{.3\textwidth}
\includegraphics[width=\textwidth,
keepaspectratio]{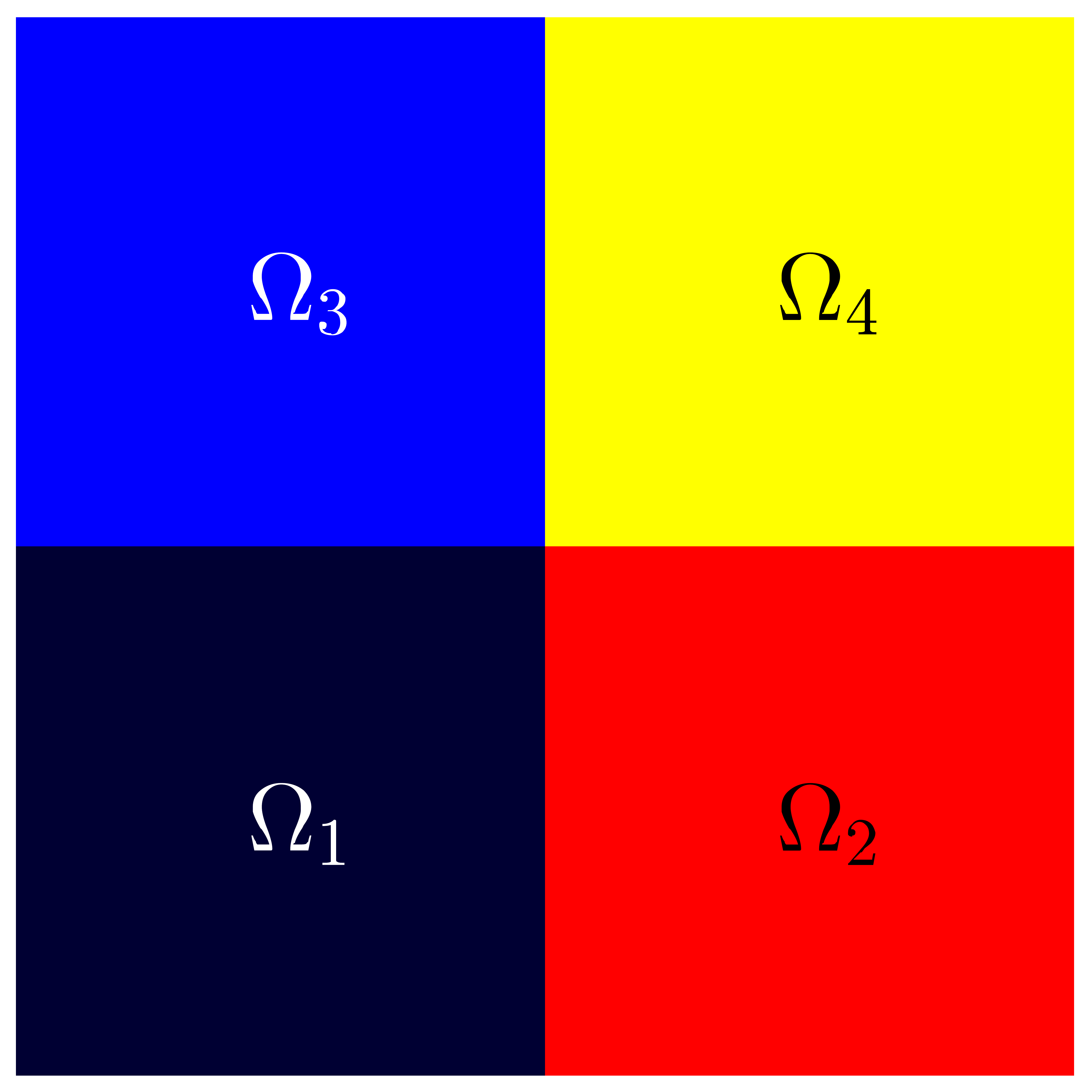}
\end{minipage}
\hspace{4mm}
\centering
\begin{minipage}[c]{.38\textwidth}
\includegraphics[width=\textwidth,
keepaspectratio]{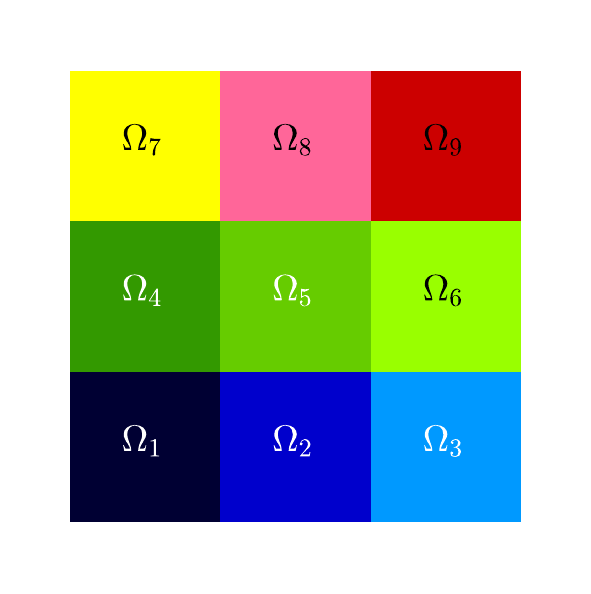}
\end{minipage}
\caption{The considered decompositions of the spatial domain $D = [0,1]^2$. \label{fig:problem_domains}}
\end{figure}

In this section we describe the test case we considered to asses the numerical performance of the diverse algorithms. Let $D = [0,1]^2 =\cup_{k=1}^K \Omega_k$ be a decomposition of the spatial domain. We consider problem \eqref{eq:stoc_ell_pde} with $f \equiv 1$ and 
$$
a(x;y) = \sum_{k=1}^K y_k \mathbbm{1}_{\Omega_k}(x), \qquad \textrm{for $x\in D$}
$$
and $y=(y_1,\dots,y_K)\in \Gamma = [1,3]^K$.
In other words, we are considering a diffusion problem on $D$, where the diffusion coefficient is constant on the elements of a partition of $D$ in $K$ zones.
We study the case of a stochastic parameter $\mathbf{Y}=(Y_1,\dots,Y_K)$ where $Y_k$'s are independent random variables with a shifted and re-scaled Beta distribution:
$$
Y_k \sim 2\cdot\mathrm{Beta}(\alpha_k,\beta_k) + 1,
$$
for some positive distribution parameters $\alpha_k,\beta_k$, $k=1,\dots,K$. We consider a uniform decomposition of the domain $D$ for either $K=4$ (low-dimensional case) or $K=9$ (higher-dimensional case), as illustrated in Figure \ref{fig:problem_domains}.

\subsection{Weighted reduced order methods}

We implemented the standard (non-weighted) and the weighted greedy algorithm for construction of reduced order spaces. We took $\omega = \sqrt{\rho}$ in the weighted case (this choice being motivated by \eqref{eq:w_gre_sqrt_est}) and we sampled $\Xi_t$ using diverse techniques: 
\begin{itemize}
\item sampling from uniform distribution;
\item sampling from equispaced grid;
\item sampling from the (Beta) distribution.
\end{itemize}
We choose $y^1 = (2,\dots,2)\in\mathbb{R}^K$ as first parameter in Algorithm \ref{alg:w_gre} (the mode of the distribution $\mathbf{Y}$). The size of the training set was chosen to be $n_t = 1000$ for the case $K = 4$ and $n_t = 2000$ for the case $K = 9$.

Furthermore, we implemented diverse versions of the weighted POD algorithm and we compared their performance. As explained above, different version of the weighted POD algorithm are based on different quadrature formulae. Each formula is specified by a set of nodes $\Xi_t = \bra{y_1,\dots,y_{n_t}}$ and a respective set of weights $w_i$ (note that these can be different from the weights $\omega_i$ of the quadrature formula; below we also denote $\rho_i=\rho(y_i)$). In particular, we experimented the following weighted variants of the POD algorithm with non-sparse grids:
\begin{itemize}
\item Uniform Monte-Carlo: the nodes are uniformly sampled and the weights are given by $w_i = \rho_i/n_t$;
\item Monte-Carlo: the nodes are sampled from the distribution of $\mathbf{Y}$ and uniformly weighted;
\item Clenshaw-Curtis (Gauss-Legendre, respectively): the nodes are the ones of  Clen- shaw-Curtis (Gauss-Legendre, resp.) tensor product quadrature formula and the weights are given by $w_i = \rho_i\omega_i$;
\item Gauss-Jacobi: the nodes are the ones of the Gauss-Jacobi tensor product quadrature formula and the weights are given by $w_i = \omega_i$.
\end{itemize}

In the low-dimensional case ($K=4$) we tested the methods for two different values of the training set size, $n^1_t$ and $n^2_t$. In the first test we chose $\Xi_t=\Xi_t^1$ to be the smallest\footnote{When using tensor product quadrature rule, we can not impose the cardinality of $\Xi_t$ a priori. We also note that when we use the Clenshaw-Curtis approximation, the majority of the points in $\Xi_t$ lies on $\partial\Gamma$: these point are completely negligible, since $\rho|_{\partial\Gamma} \equiv 0$. So, in this case, we need to take a considerably larger value for $n_t$ to reach the desired cardinality of nodes in the interior $\mathring{\Gamma}$.} possible set such that $\abs{\Xi_t^1 \setminus \partial\Gamma}\geq 100$ and in the second test we chose $\Xi_t=\Xi_t^2$ to be the smallest possible set such that $\abs{\Xi_t^2 \setminus \partial\Gamma}\geq 500$. For the higher-dimensional case ($K=9$), the sizes of the sets of nodes of Gauss-Legendre and Gauss-Jacobi rules was $n_t=2^9=512$; indeed, the next possible choice $n_t=3^9=19683$ was computationally impracticable. Moreover, Clenshaw-Curtis formula was not test, because a very large training set would have been required as $\abs{\Xi_t \cap \mathring{\Gamma}}=1$ for (the already impractible choice) $n_t=3^9$. In this case, the use of sparse Gauss-Legendre/Gauss-Jacobi quadrature rules provided a more representative set of nodes consisting of just $n_t = 181$ nodes. 
A summary of the formulae and training set sizes used is reported in Table \ref{table:4dz_pod}.
\begin{table}[tbp]
\centering
{\renewcommand\arraystretch{1.2}
\begin{tabular}{|l|c|c|c|c|}
\hline
 & \multirow{2}{*}{$w_i$} & \multicolumn{2}{|c|}{$K=4$} & \multicolumn{1}{|c|}{$K=9$} \\
\cline{3-5}
 & & $n_t^1$  & $n_t^2$ & $n_t$ \\
\hline
Standard & $1/n_t$ & $100$ & $500$ & $500$ \\
\hline
Monte-Carlo & $1/n_t$ & $100$ & $500$ & $500$ \\
\hline
Uniform Monte-Carlo & $\rho_i/n_t$ & $100$ & $500$ & $2000$ \\
\hline
Clenshaw-Curtis & $\omega_i\rho_i$ & $1296$ & $2401$ & $-$\\
\hline
Gauss-Legendre & $\omega_i\rho_i$ & $256$ & $625$ & $512$ \\
\hline
Gauss-Jacobi & $\omega_i$ & $256$ & $625$ & $512$ \\
\hline
Sparse Gauss-Legendre & $\omega_i\rho_i$ & $-$ & $-$ & $181$ \\
\hline
Sparse Gauss-Jacobi & $\omega_i$ & $-$ & $-$ & $181$ \\
\hline
\end{tabular}} 
\caption{Weights and sizes of training sets in the diverse weighted POD variants tested.\label{table:4dz_pod}}
\end{table}

\subsection{Results}

We compared the performance of the weighted greedy and weighted POD algorithms for computing the expectation of the solution. In particular, we computed the error \eqref{eq:w_pod_min_qua}
using a Monte Carlo approximation, i.e.,
\begin{equation}
\label{eq:w_pod_err}
\norm{u(\mathbf{Y})-u_N(\mathbf{Y})}^2_\mathbb{H}
\simeq \frac{1}{M} \sum_{m=1}^M \norm{u(y^m)-u_N(y^m)}^2_\mathbb{V},
\end{equation}
where $y^1,\dots,y^M$ are $M$ independent realizations of $\mathbf{Y}$. For truth problem solutions, we adopted $\mathbb{V}_\delta$ to be the classical $\mathbb{P}^1$-FE approximation space. Three cases were carried out:
\begin{enumerate}
\item distribution parameters $(\alpha_i,\beta_i)=(10,10)$, for $i = 1, \hdots, K = 4$;
\item distribution parameters $(\alpha_i,\beta_i)=(10,10)$, for $i = 1, \hdots, K = 9$;
\item distribution parameters $(\alpha_i,\beta_i)=(75,75)$, for $i = 1, \hdots, K = 9$, resulting in a more concentrated distribution than case 2.
\end{enumerate}

\begin{figure}[p]
\centering
\begin{minipage}[c]{0.45\textwidth}
\includegraphics[width=\textwidth,
keepaspectratio]{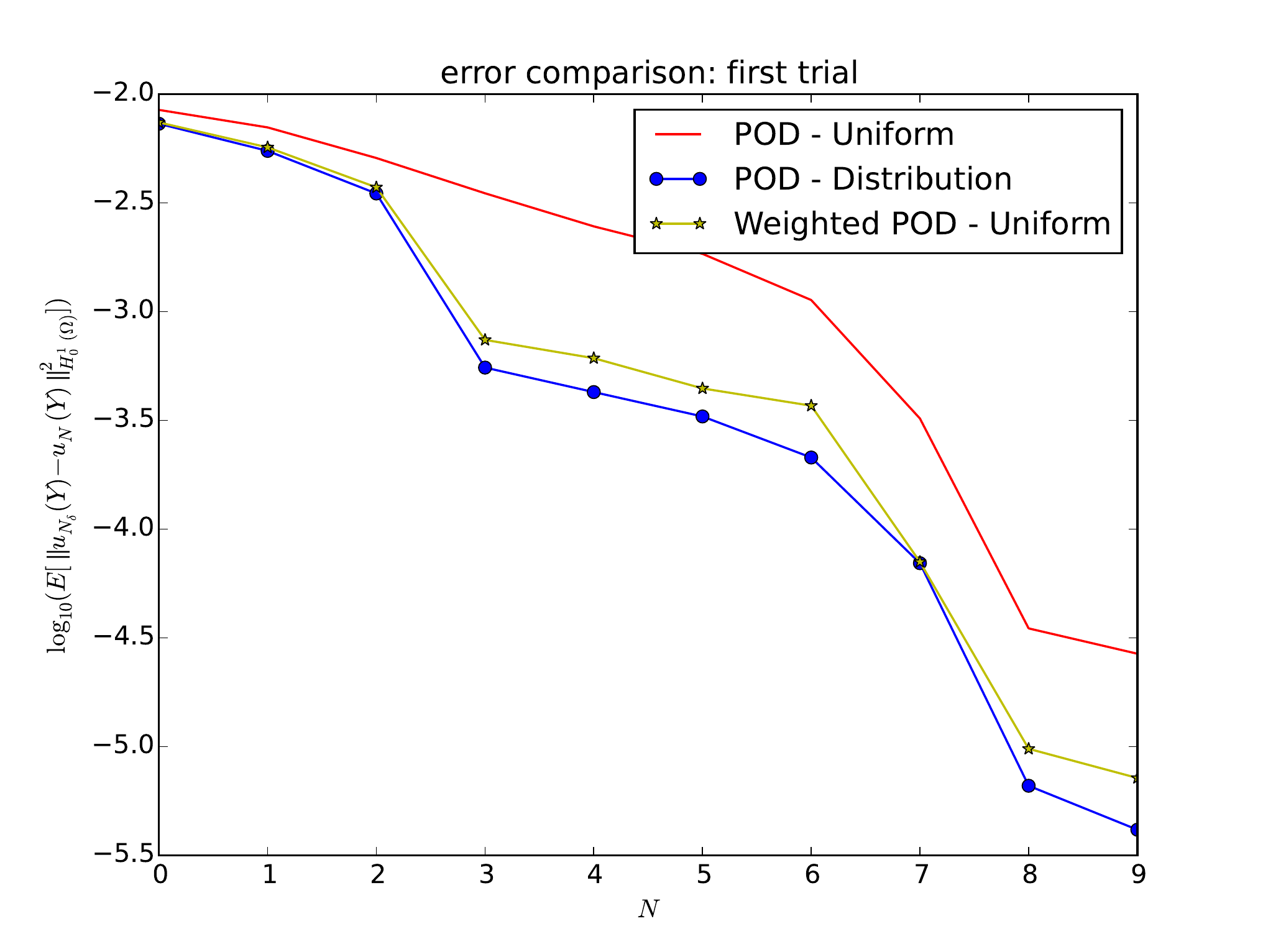}
\end{minipage}
\hspace{4mm}
\begin{minipage}[c]{.45\textwidth}
\includegraphics[width=\textwidth,
keepaspectratio]{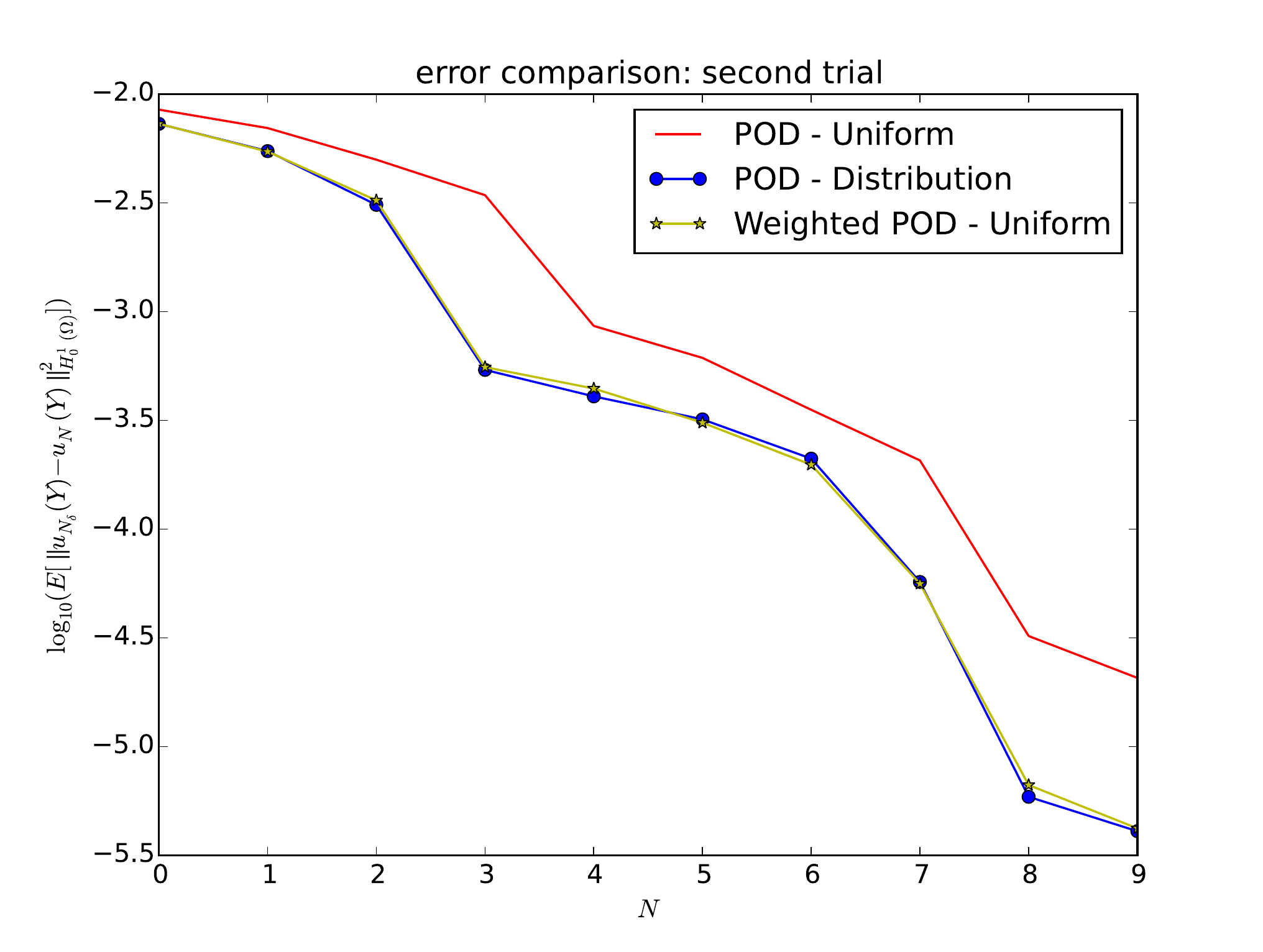}
\end{minipage}
\caption{\emph{Low-dimensional case ($K=4$).} Plots of the error \eqref{eq:w_pod_err} obtained using standard (POD - uniform), uniform Monte-Carlo (Weighted POD - Uniform) and Monte-Carlo (POD - distribution) POD algorithms. Left: $\Xi_t^1$. Right: $\Xi_t^2$.\label{figure:LOW:wPOD}}
\end{figure}

\begin{figure}[p]
\centering
\begin{minipage}[c]{.45\textwidth}
\includegraphics[width=\textwidth,
keepaspectratio]{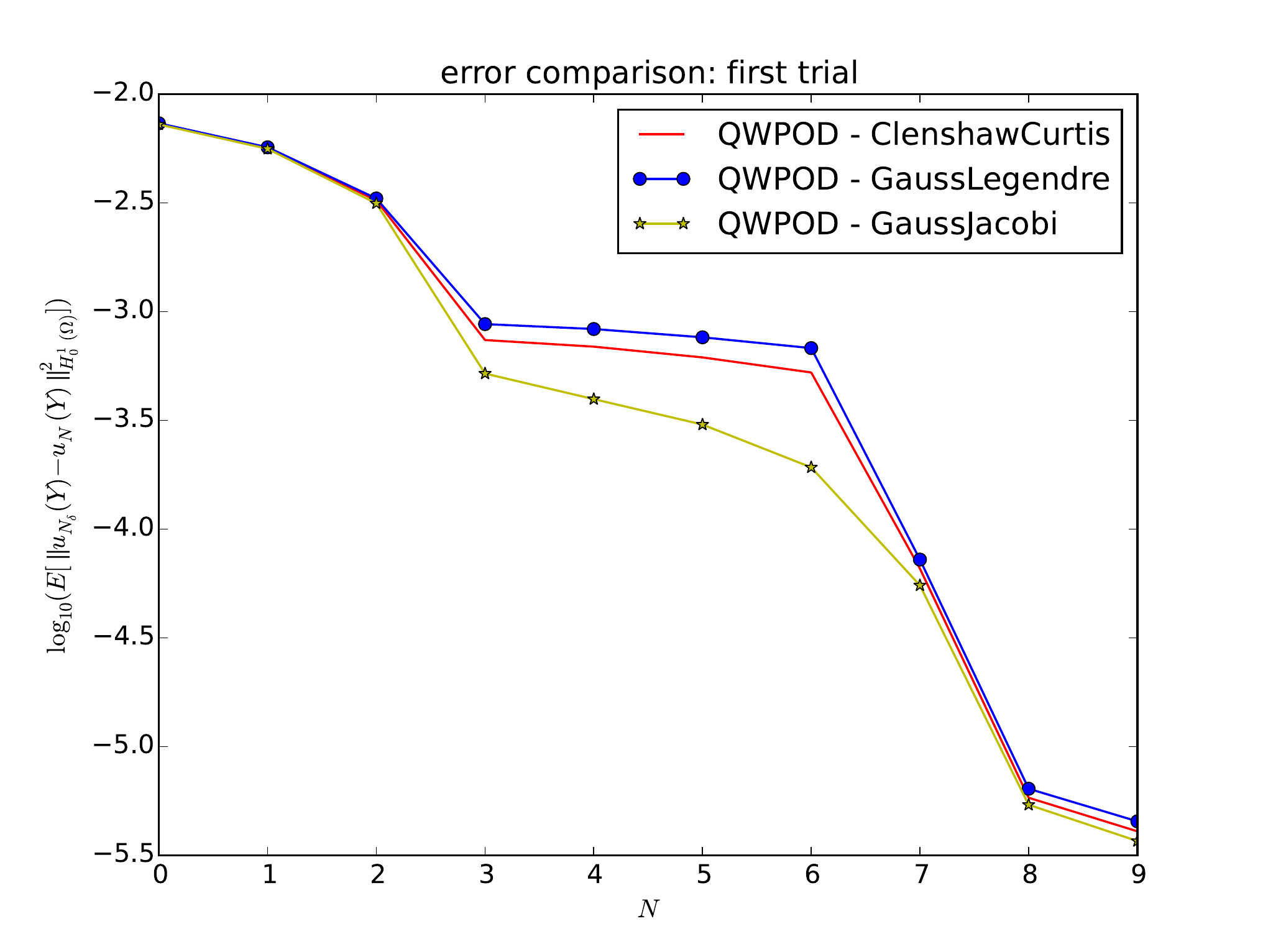}
\end{minipage}
\hspace{4mm}
\begin{minipage}[c]{.45\textwidth}
\includegraphics[width=\textwidth,
keepaspectratio]{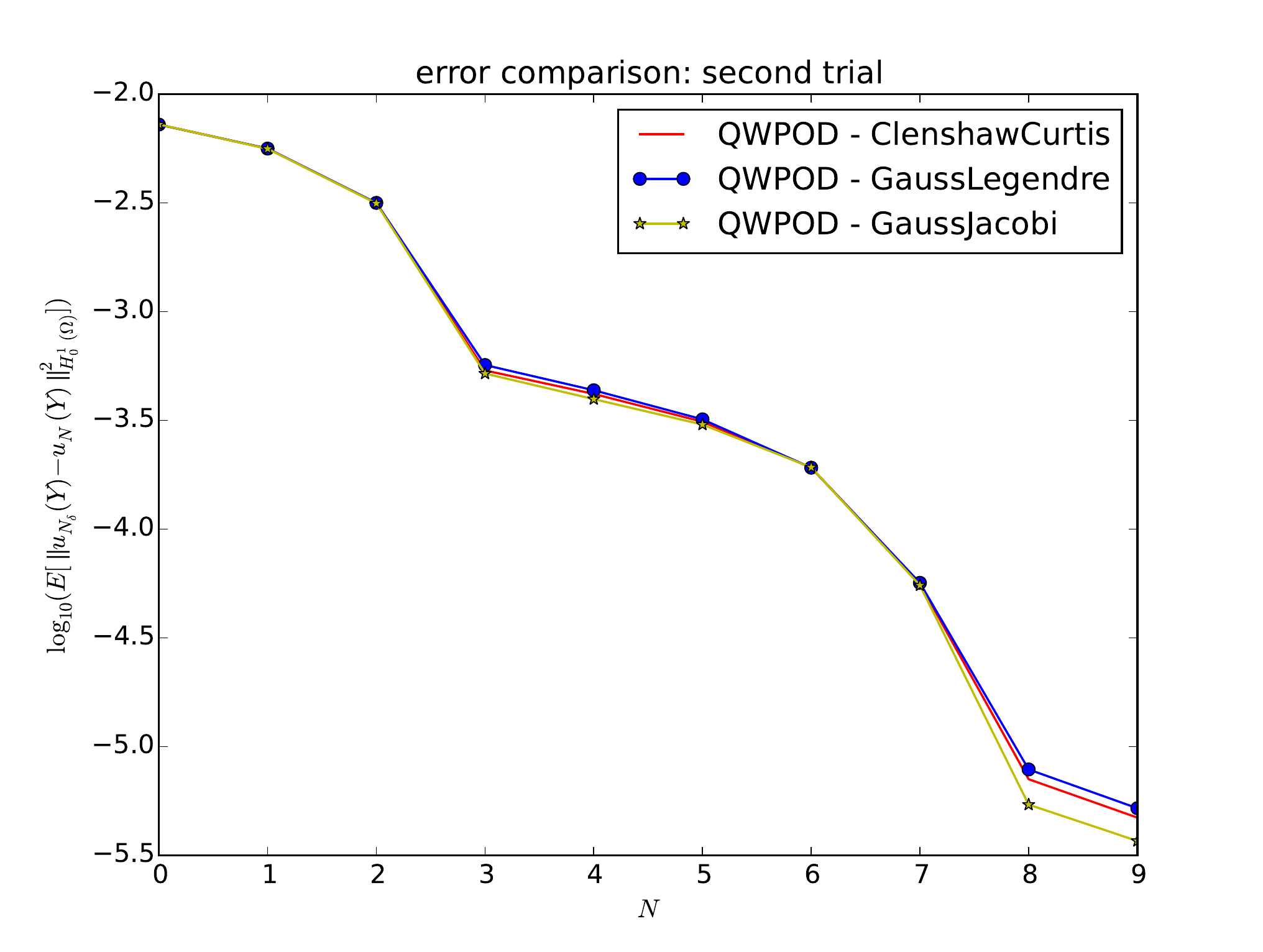}
\end{minipage}
\caption{\emph{Low-dimensional case ($K=4$).} Plots of the error \eqref{eq:w_pod_err} obtained using Clenshaw-Curtis, Gauss-Legendre and Gauss-Jacob POD algorithms. Left: $\Xi_t^1$. Right: $\Xi_t^2$.\label{figure:LOW:qPOD}}
\end{figure}

\begin{figure}[p]
\centering
\begin{minipage}[c]{0.45\textwidth}
\includegraphics[width=\textwidth,
keepaspectratio]{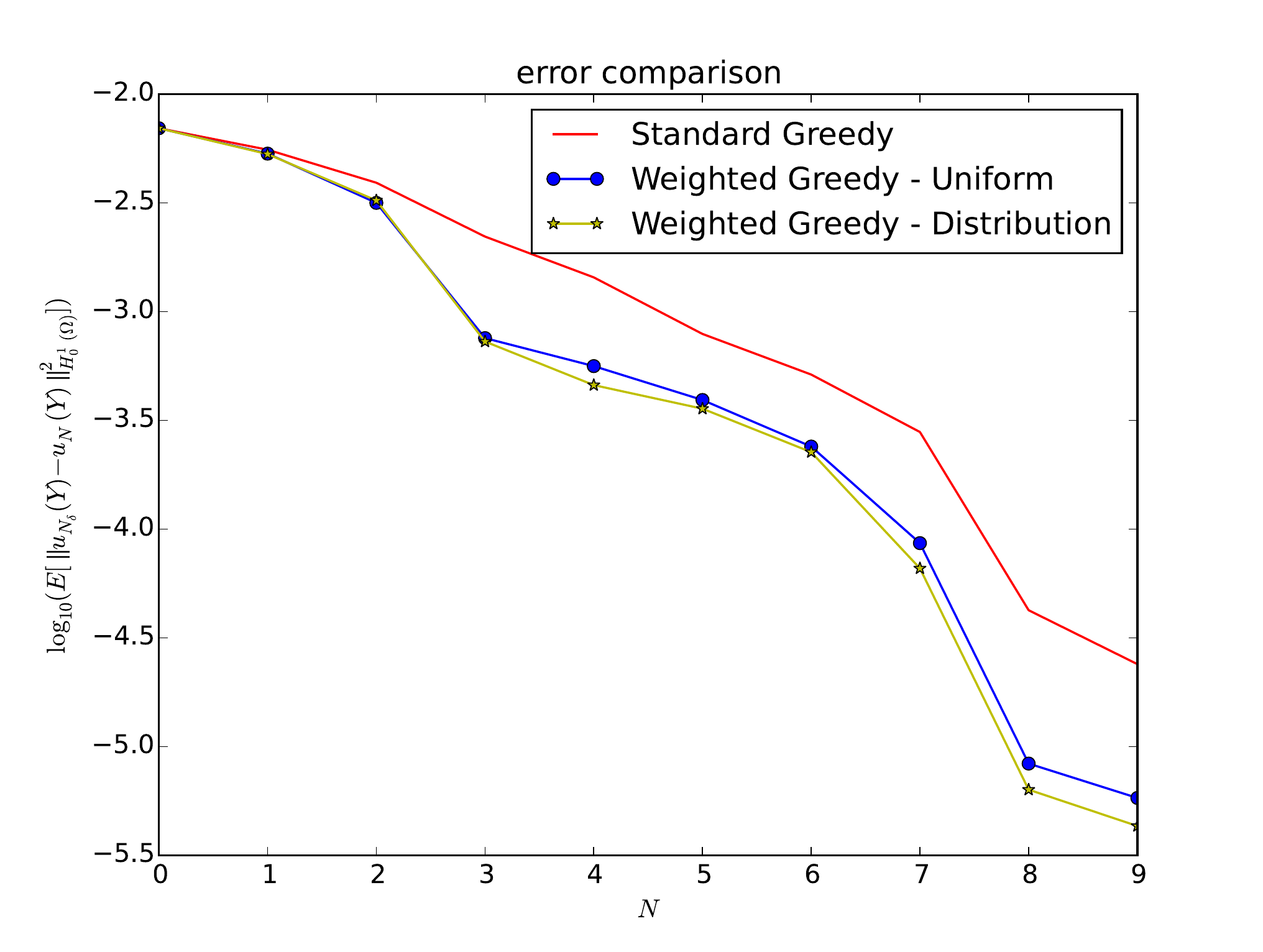}
\end{minipage}
\hspace{4mm}
\begin{minipage}[c]{0.45\textwidth}
\includegraphics[width=\textwidth,
keepaspectratio]{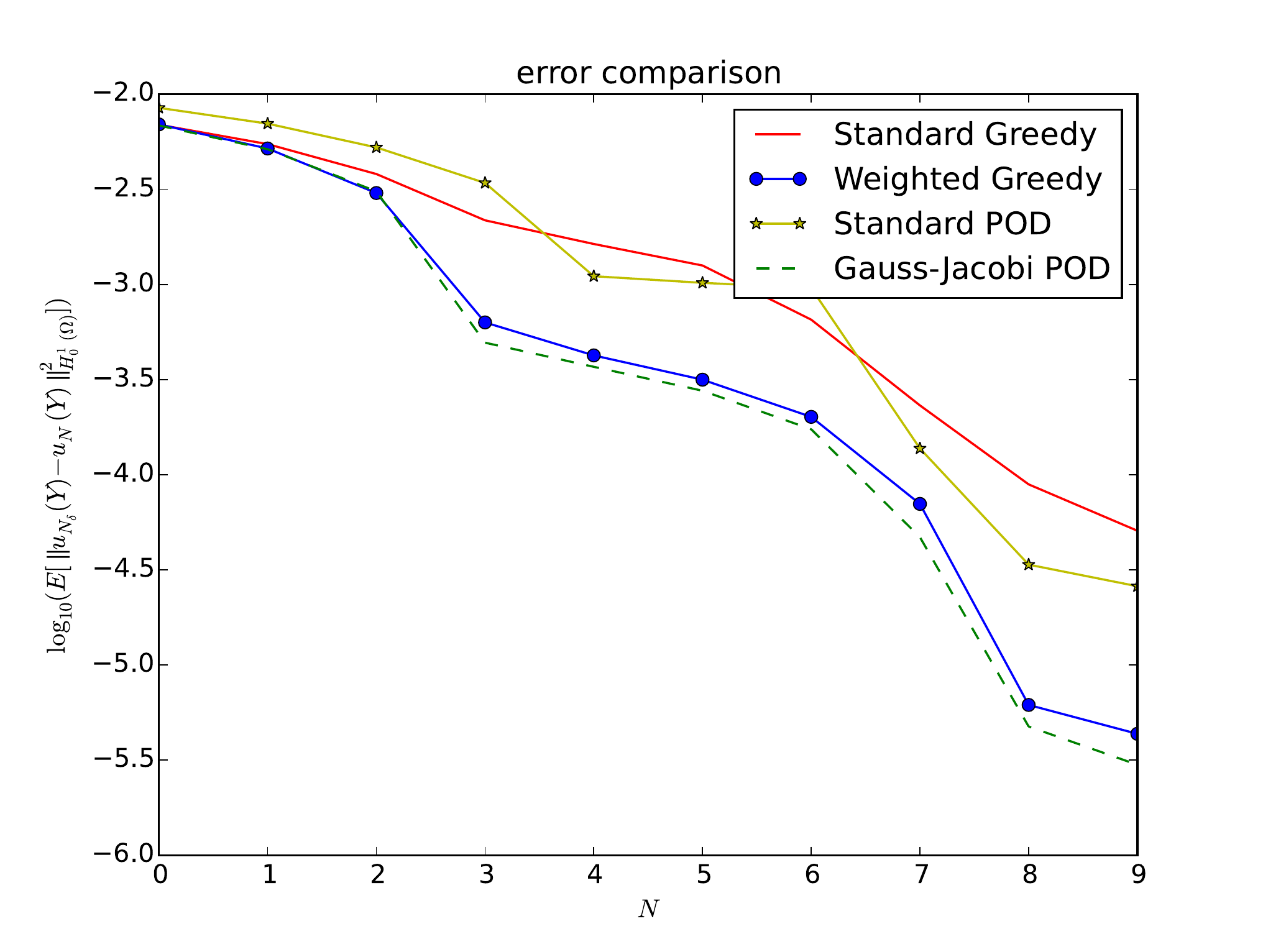}
\end{minipage}
\caption{\emph{Low-dimensional case ($K=4$).} Left: Plots of the error \eqref{eq:w_pod_err} obtained using standard greedy algorithm and weighted greedy algorithm with uniform sampling and sampling of the distribution. Right: Plots of the error \eqref{eq:w_pod_err} obtained using standard and weighted greedy and standard and Gauss-Jacobi POD algorithms. \label{figure:LOW:WG+MIX}}
\end{figure}

\begin{figure}[p]
\centering
\begin{minipage}[c]{.45\textwidth}
\includegraphics[width=\textwidth,
keepaspectratio]{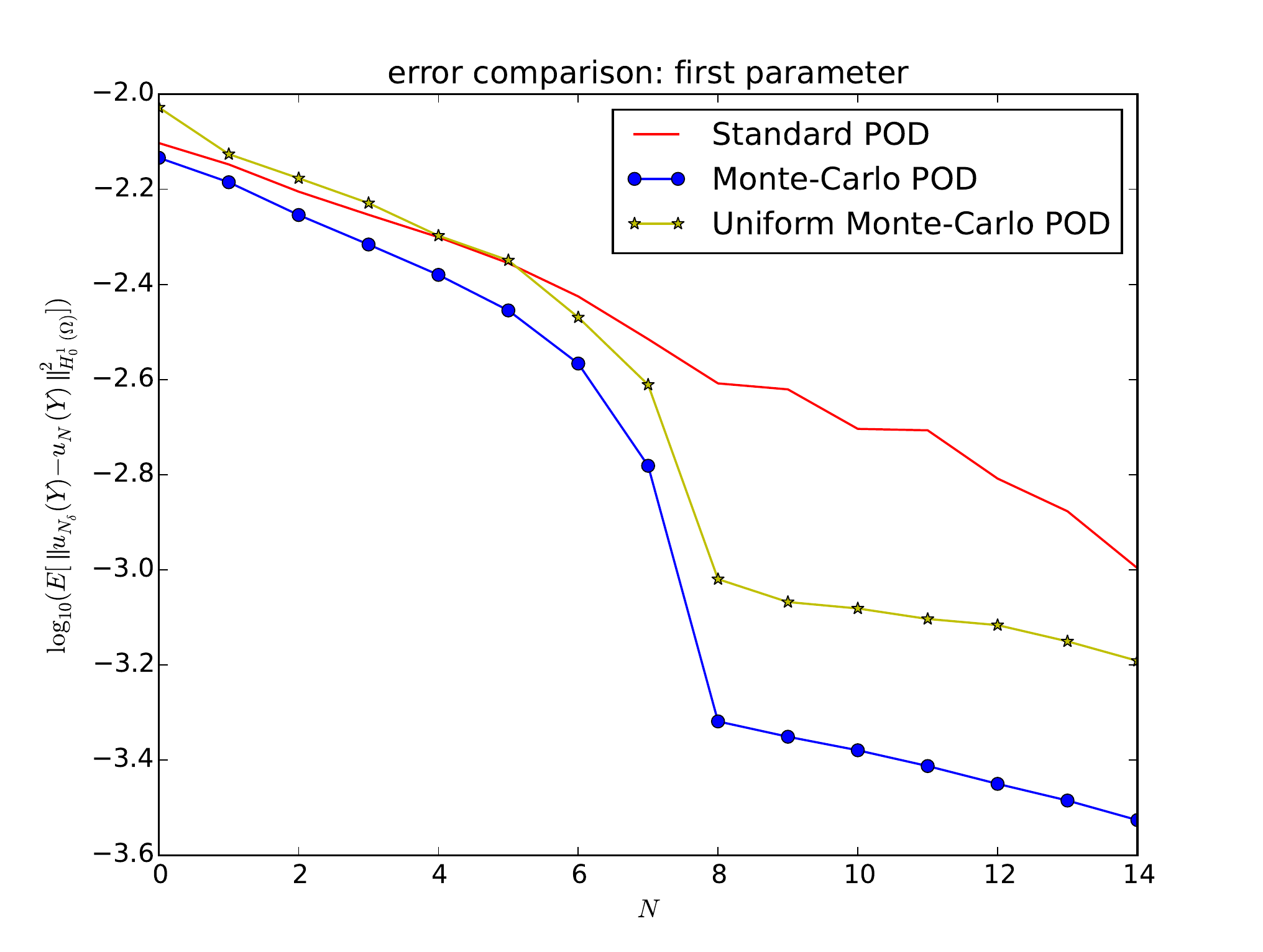}
\end{minipage}
\hspace{4mm}
\begin{minipage}[c]{.45\textwidth}
\includegraphics[width=\textwidth,
keepaspectratio]{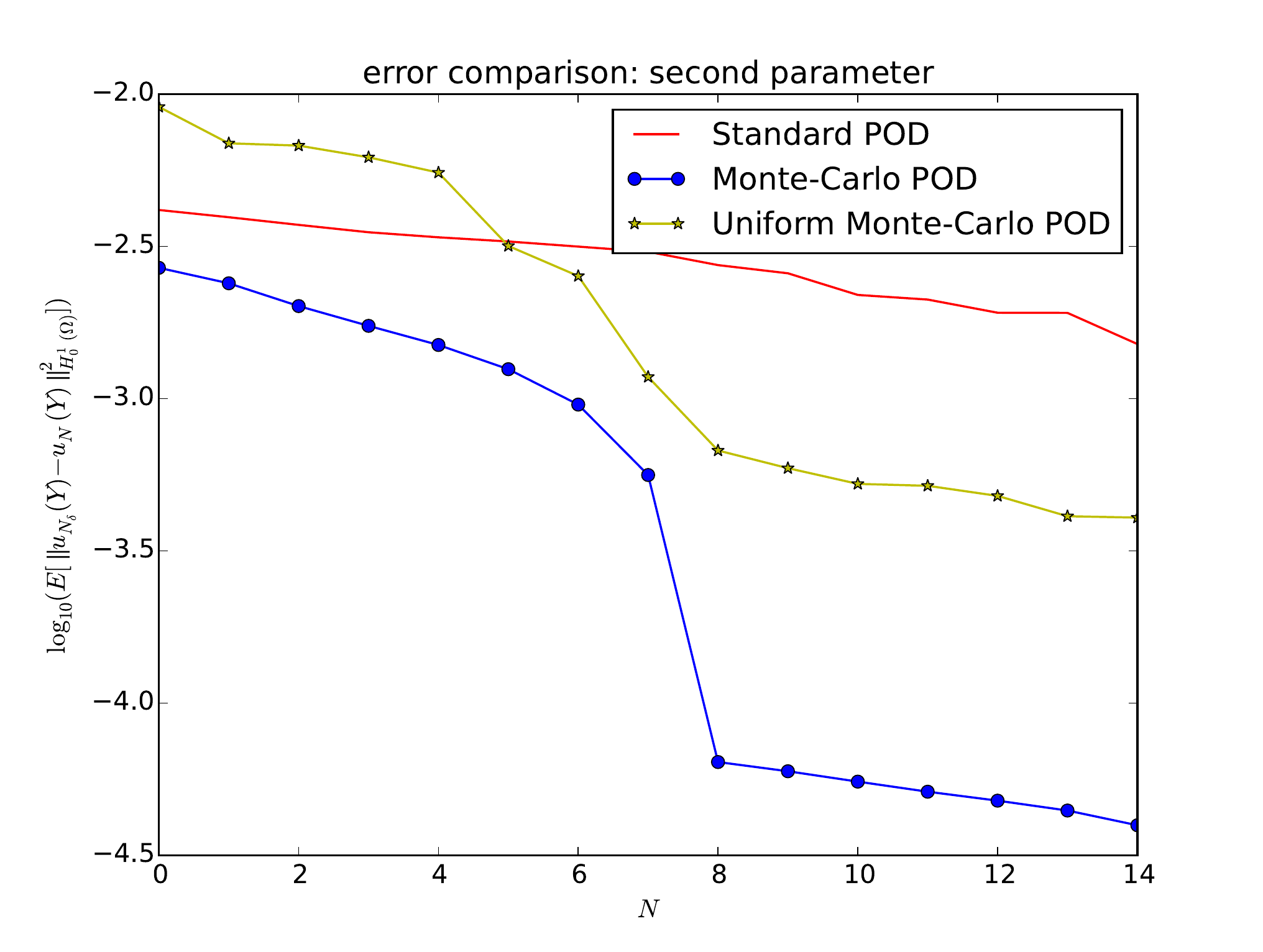}
\end{minipage}
\caption{\emph{Higher-dimensional case ($K=9$).} Plots of the error \eqref{eq:w_pod_err} obtained using standard, uniform Monte-Carlo and Monte-Carlo POD algorithms. Left: $\alpha_i=\beta_i=10$. Right: $\alpha_i=\beta_i=75$. \label{figure:HIGH:wPOD}}
\end{figure}

\begin{figure}[p]
\centering
\begin{minipage}[c]{.45\textwidth}
\includegraphics[width=\textwidth,
keepaspectratio]{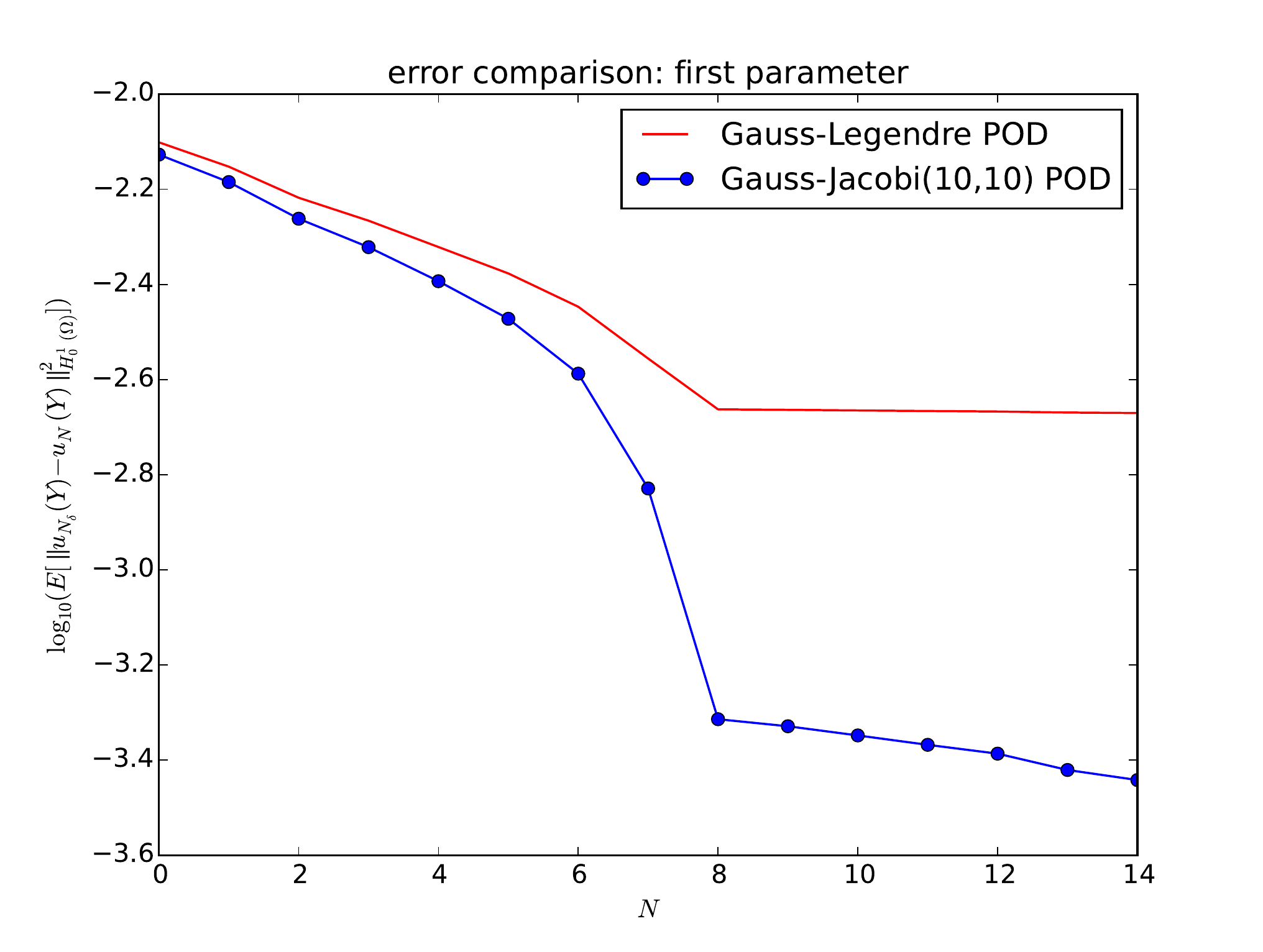}
\end{minipage}
\hspace{4mm}
\begin{minipage}[c]{.45\textwidth}
\includegraphics[width=\textwidth,
keepaspectratio]{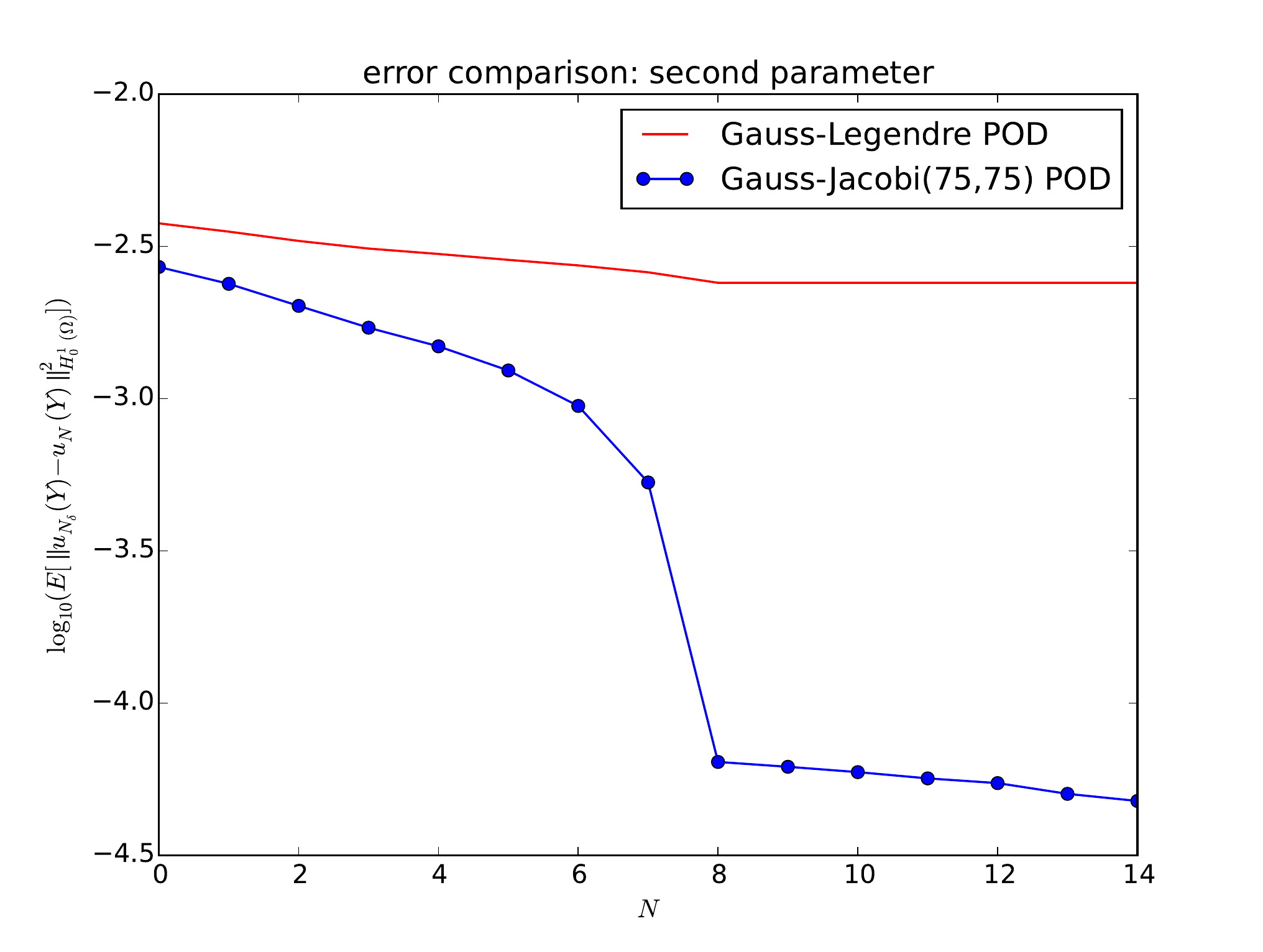}
\end{minipage}
\caption{\emph{Higher-dimensional case ($K=9$).} Plots of the error \eqref{eq:w_pod_err} obtained using (tensor product) Gauss-Legendre and Gauss-Jacobi POD algorithms. Left: $\alpha_i=\beta_i=10$. Right: $\alpha_i=\beta_i=75$. \label{figure:HIGH:qPOD}}
\end{figure}

\begin{figure}[p]
\centering
\begin{minipage}[c]{.45\textwidth}
\includegraphics[width=\textwidth,
keepaspectratio]{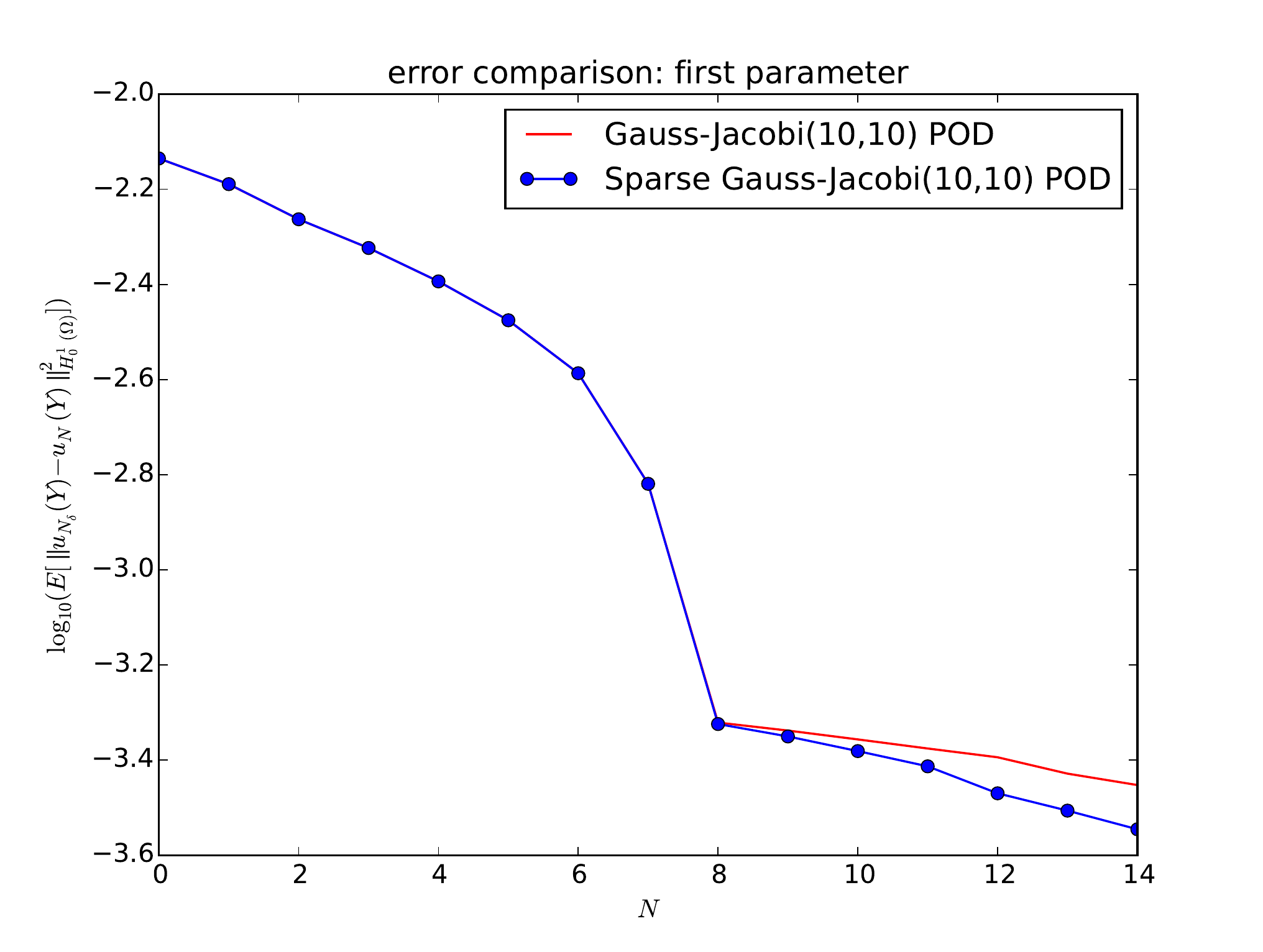}
\end{minipage}
\hspace{4mm}
\begin{minipage}[c]{.45\textwidth}
\includegraphics[width=\textwidth,
keepaspectratio]{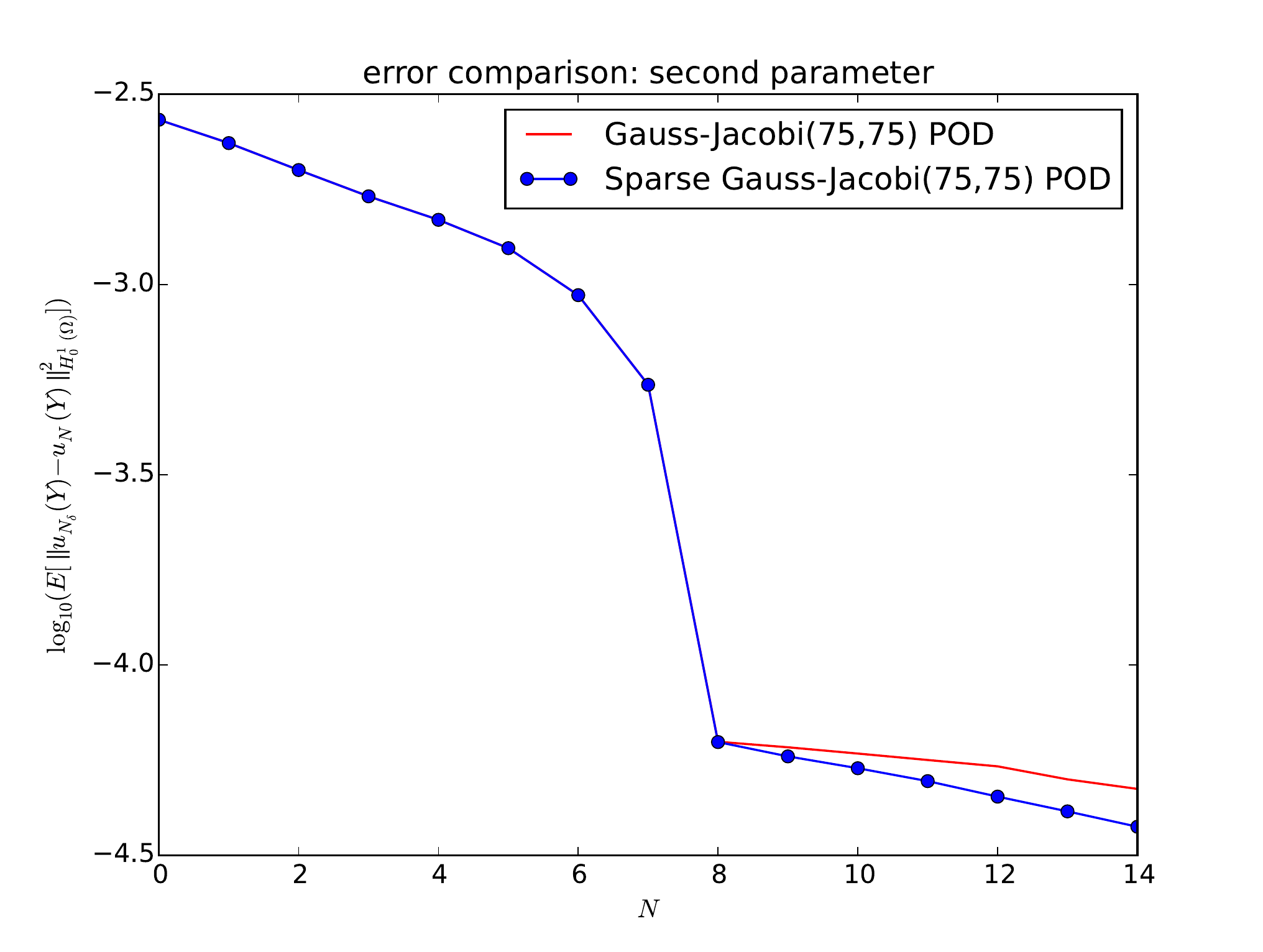}
\end{minipage}
\caption{\emph{Higher-dimensional case ($K=9$).} Plots of the error \eqref{eq:w_pod_err} obtained using Gauss-Legendre POD and sparse Gauss-Legendre POD algorithms. Left: $\alpha_i=\beta_i=10$. Right: $\alpha_i=\beta_i=75$. \label{figure:HIGH:sPOD}}
\end{figure}

\begin{figure}[p]
\centering
\begin{minipage}[c]{.45\textwidth}
\includegraphics[width=\textwidth,
keepaspectratio]{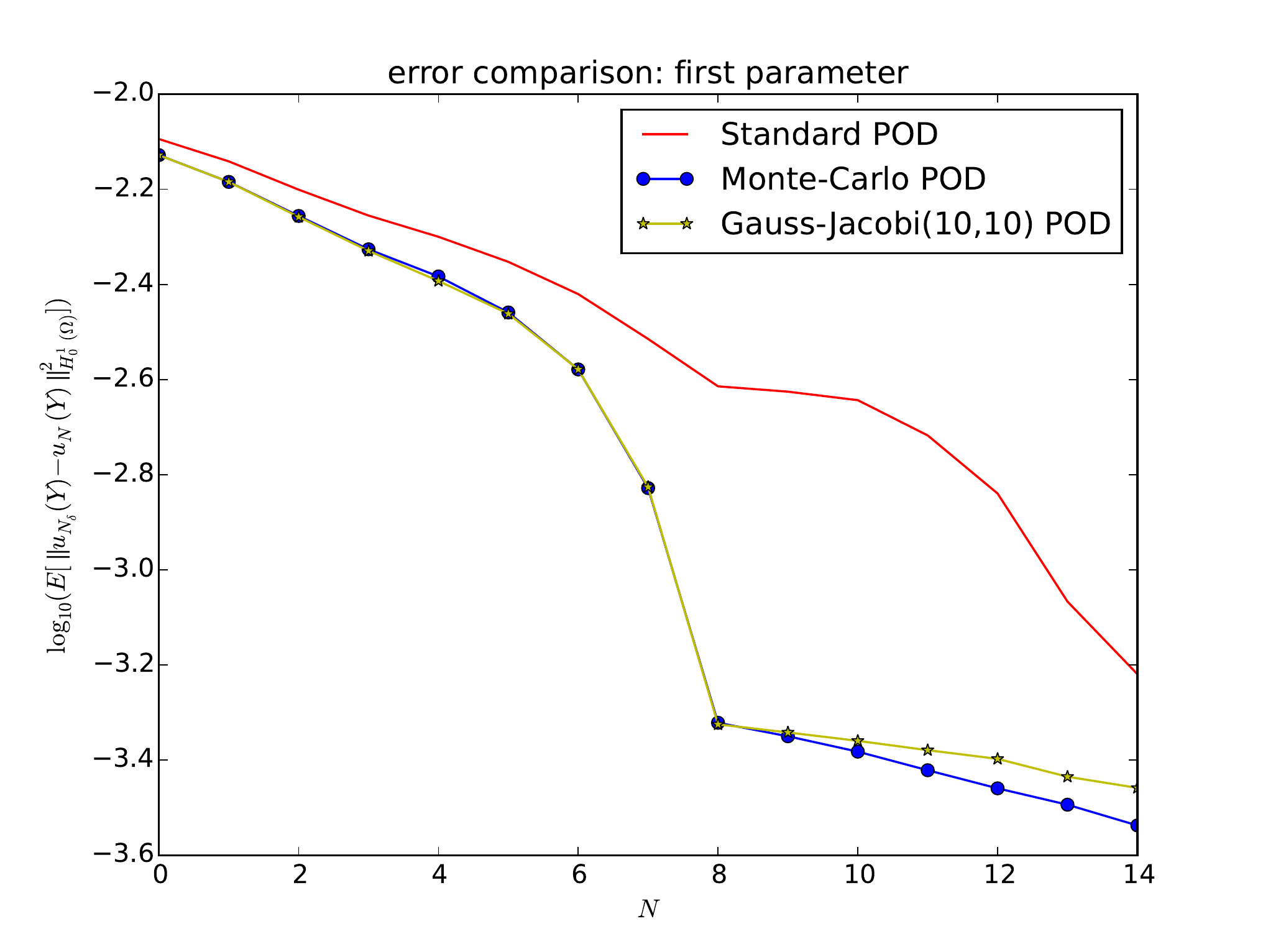}
\end{minipage}
\hspace{4mm}
\begin{minipage}[c]{.45\textwidth}
\includegraphics[width=\textwidth,
keepaspectratio]{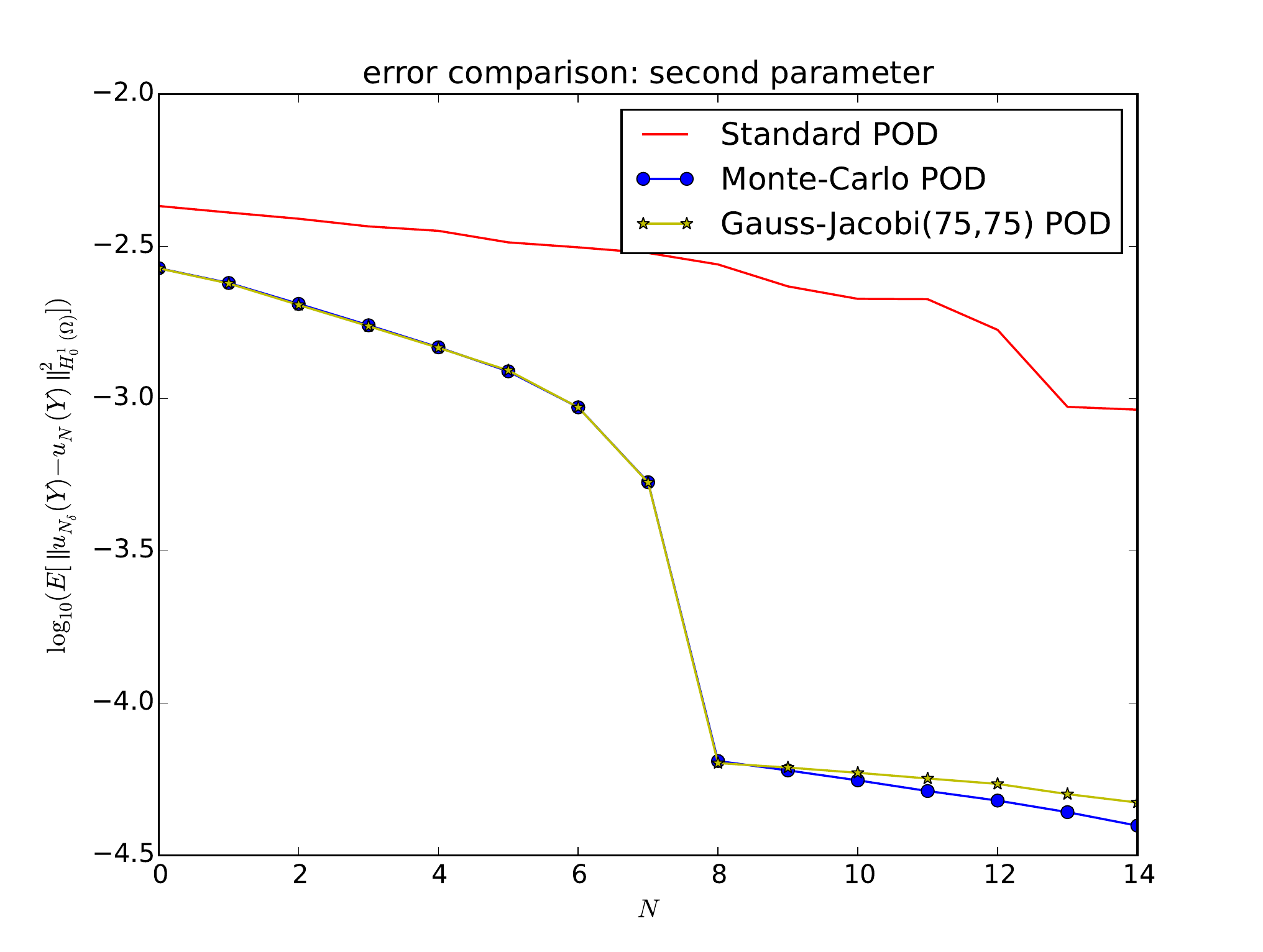}
\end{minipage}
\caption{\emph{Higher-dimensional case ($K=9$).} Plots of the error \eqref{eq:w_pod_err} obtained using standard, Monte-Carlo and Gauss-Jacobi POD algorithms. Left: $\alpha_i=\beta_i=10$. Right: $\alpha_i=\beta_i=75$. \label{figure:HIGH:cPOD}}
\end{figure}

\begin{figure}[p]
\centering
\begin{minipage}[c]{.45\textwidth}
\includegraphics[width=\textwidth,
keepaspectratio]{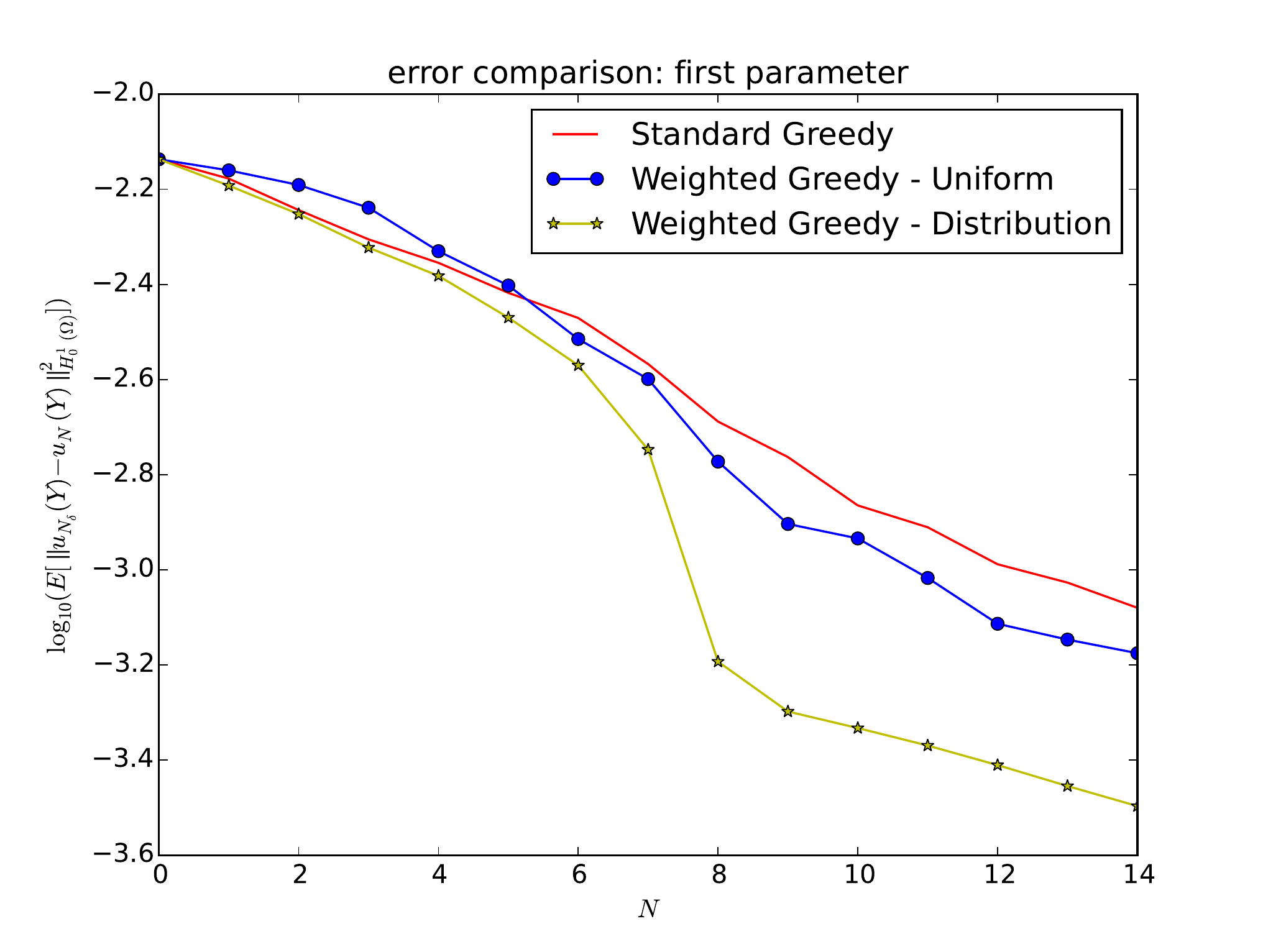}
\end{minipage}
\hspace{4mm}
\begin{minipage}[c]{.45\textwidth}
\includegraphics[width=\textwidth,
keepaspectratio]{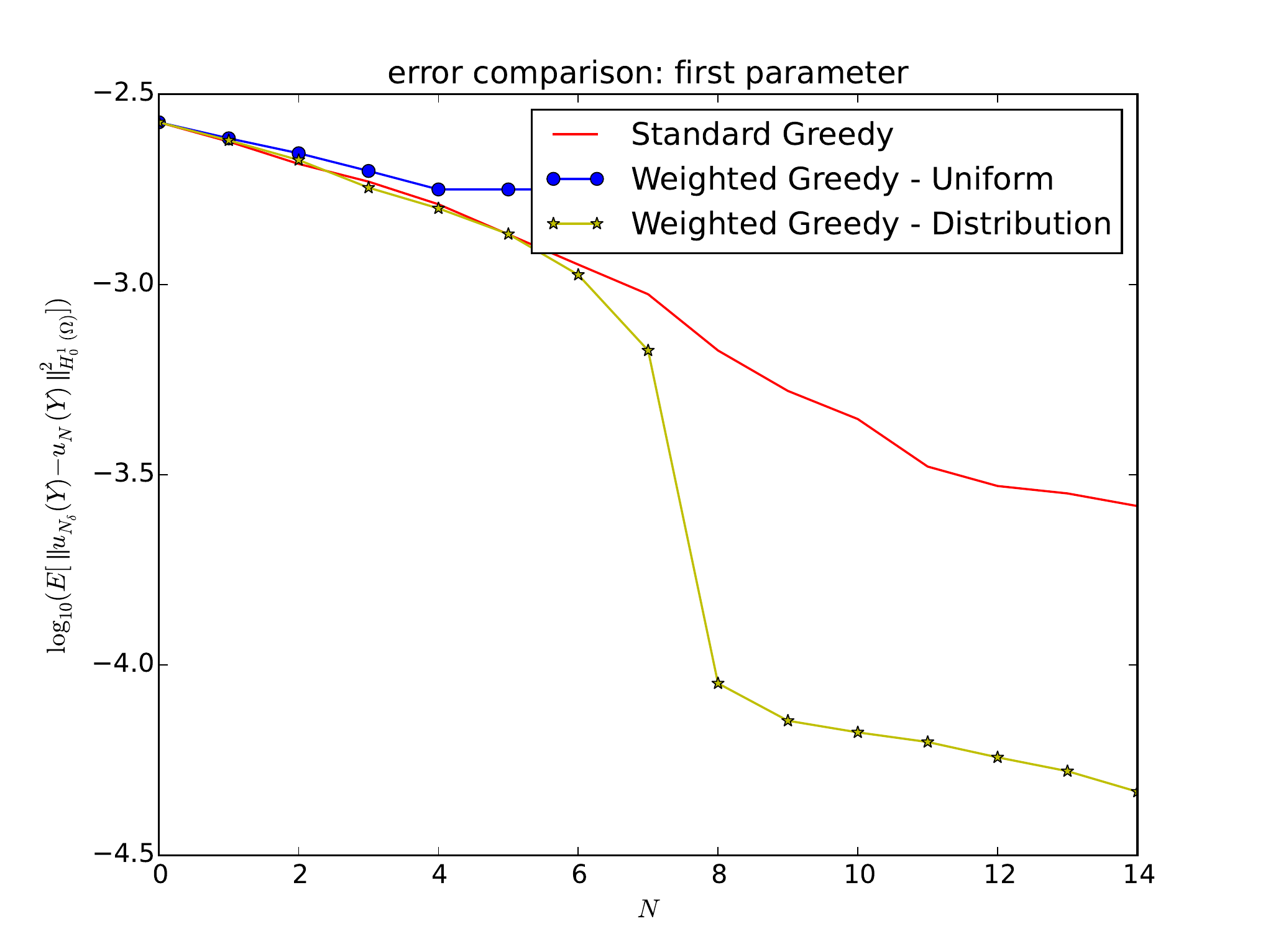}
\end{minipage}
\caption{\emph{Higher-dimensional case ($K=9$).} Plots of the error \eqref{eq:w_pod_err} obtained using standard and weighted greedy algorithms with uniform sampling and from the distribution. Left: $\alpha_i=\beta_i=10$. Right: $\alpha_i=\beta_i=75$. \label{figure:HIGH:GR}}
\end{figure}

\begin{figure}[p]
\centering
\begin{minipage}[c]{.45\textwidth}
\includegraphics[width=\textwidth,
keepaspectratio]{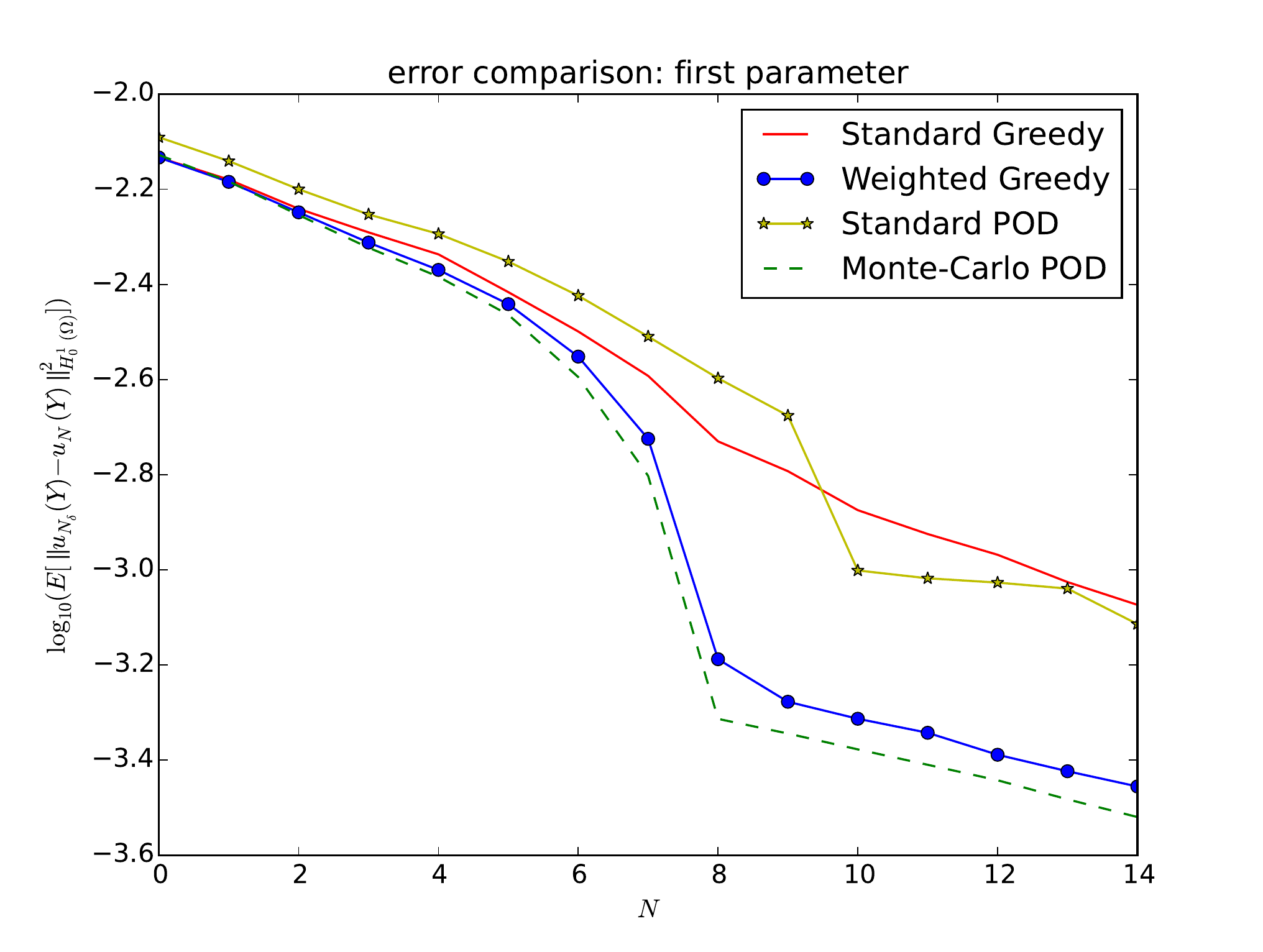}
\end{minipage}
\hspace{4mm}
\begin{minipage}[c]{.45\textwidth}
\includegraphics[width=\textwidth,
keepaspectratio]{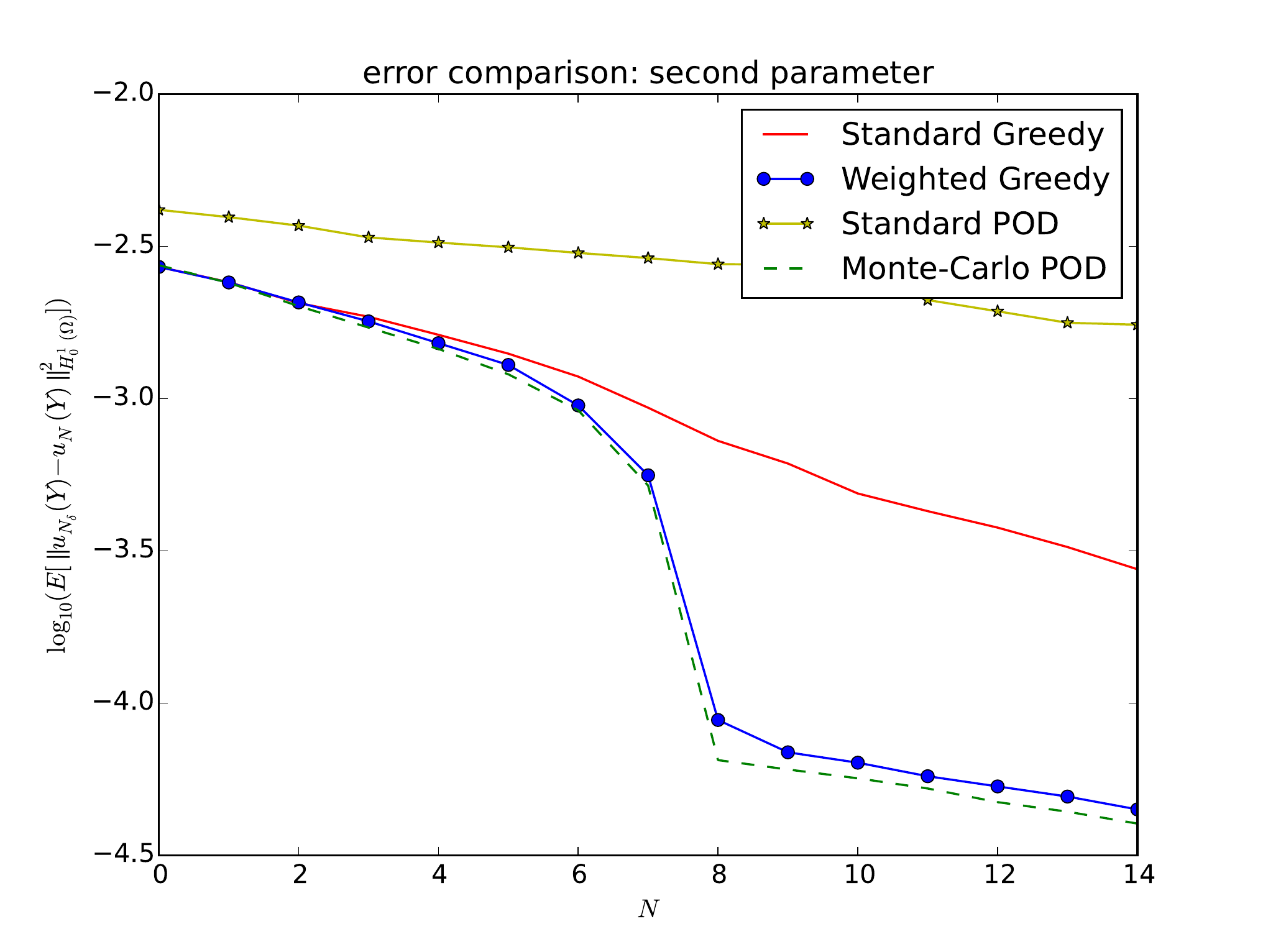}
\end{minipage}
\caption{\emph{Higher-dimensional case ($K=9$).} Plots of the error \eqref{eq:w_pod_err} obtained using standard and weighted Greedy and standard and Monte-Carlo POD algorithms. Left: $\alpha_i=\beta_i=10$. Right: $\alpha_i=\beta_i=75$. \label{figure:HIGH:MIX}}
\end{figure}

Figures \ref{figure:LOW:wPOD}--\ref{figure:HIGH:MIX} show the graphs of the error \eqref{eq:w_pod_err} (in a $\log_{10}$ scale), as a function of the reduced order space dimension $N$, for different methods. In particular, Figs 3--5 collect the results for case 1, as well as Figs. 6--8 those of case 2, while results for case 3 are shown in Figs 9--11. From these plots, we can gather the following conclusions:
\begin{itemize}
\item \emph{Weighted vs standard}. The weighted versions of the POD and greedy algorithms outperform the performance of their standard counter-parts in the stochastic setting, see Figures \ref{figure:LOW:wPOD} (POD) and \ref{figure:LOW:WG+MIX} (Greedy) for case 1. Such a difference is even more evident in for higher parametric dimension (compare Figure \ref{figure:LOW:wPOD} to Figure \ref{figure:HIGH:wPOD} for the POD case) or in presence of higher concentrated parameters distributions (compare Figure \ref{figure:HIGH:wPOD} to Figure \ref{figure:HIGH:cPOD} for the POD case, and see Figures \ref{figure:HIGH:GR} and \ref{figure:HIGH:MIX} for the greedy case).
\item \emph{Importance of representative training set}.
The Monte-Carlo and Gauss-Jacobi POD algorithms outperform the other weighted POD variants (Figures \ref{figure:LOW:qPOD} for case 1, as well as Figures \ref{figure:HIGH:wPOD} and \ref{figure:HIGH:qPOD} for case 2). In the low-dimensional case 1, we can still recover the same accuracy also with the other weighted variants, at the cost of using a larger training set $\Xi_t$ (Figure \ref{figure:LOW:qPOD}). However, for the higher dimensional case 2, the choice of the nodes plays a much more fundamental role (see Figures \ref{figure:HIGH:wPOD} and \ref{figure:HIGH:qPOD}). Monte-Carlo and Gauss-Jacobi methods perform significantly better because the underlying quadrature rule is designed for the specific distribution $\rho(y)\,dy$ (a Beta distribution). For more concentrated distributions, methods lacking a representative training set may eventually lead to a very inaccurate sampling of the parameter space, even resulting in numerically singular reduced order matrices despite orthonormalization of the snapshots. This is because the subspace built by the reduced order methods is a subspace of $\bra{u(y) \st y \in \Xi_t }$. Thus, if $\Xi_t$ contains only points $y$ with low probability $P\bra{Y = y}$, adding a linear combination of solutions $\bra{u(y) \st y \in\Xi_t }$ to the reduced space does not increase the accuracy of the computed statistics. Instead the weighting procedure tends to neglect such solutions, resulting in linearly dependent vectors defining the reduced order subspace. This is can also be observed for the weighted greedy algorithm (Figure \ref{figure:HIGH:GR}).
\item \emph{Breaking curse of dimensionality through sparse grids}. The presented algorithms also suffer from the fact that for increasing parameter space dimensions $K$, the number of nodes, $n_t$, increases exponentially. In the weighted POD algorithm we can mitigate this problem by adopting a sparse quadrature rule. 
Figure \ref{figure:HIGH:sPOD} shows the performances of the (tensor product) Gauss-Jacobi POD versus its sparse version. Not only does the use of sparse grids make the method much more efficient, reducing significantly the size of the training set used (in our specific case, from $n_t = 512$ to $n_t=181$), but it also results in a negligibly higher accuracy. 
On the other hand, the importance of a well-representative training set is highlighted by the use of sparse grids. 
Sparse Gauss-Legendre results in numerically singular matrices at lower $N$ (than its tensor product counter-part); sparsity makes its associated training set less representative. For this reason, the plots of the error obtained with this method were not reported.  
\item \emph{Weighted POD vs weighted greedy}. The weighted POD seems to work slightly better than the weighted greedy (see Figure \ref{figure:LOW:WG+MIX} (right) for case 1 and Figure \ref{figure:HIGH:MIX} for case 3). This is because weighted POD is designed to minimize (in some sense) the quantity \eqref{eq:w_pod_min_qua}; however, the difference in terms of the error is practically negligible. The main difference in the two algorithms lies in the different training procedure. Thanks to the availability of an inexpensive error estimator, we are able to use large training sets for greedy algorithms, while still requiring a moderate computational load during the training phase. On the flip side, the availability of different techniques for the POD algorithms also allows control of the computational cost of this algorithm, which does not require the construction of an ad-hoc error estimator, making it more suitable to study problems for which such error estimation procedure is not yet available.
\end{itemize}


\afterpage{\clearpage}

\section{Conclusion and Perspectives}

In this work we developed a weighted POD method for elliptic PDEs. The algorithm is introduced alongside the previously developed weighted greedy algorithm. While the latter aims to minimize the error in the $L^\infty$ norm, the former minimizes an approximation of the mean squared error, which is of better interpretation when weighted. Differently from the greedy, the introduced algorithm does not require the availability of error estimation. It is instead based on a quadrature rule, which can be chosen accordingly to the parameter distribution. In particular, this allows to implement sparse quadrature rules to reduce the computational cost of the offline phase as well. \\
A numerical example on a thermal block problem was carried out to test the proposed reduced order method, for either a low dimensional parametric space or two higher dimensional space of parameters. For this problem, we assessed that the weighted POD method is an efficient alternative to the weighted greedy algorithm. The numerical tests also highlighted the importance of a training set which is representative of the underlying parameter distribution. In case of a representative rule, the sparse quadrature rule based algorithm showed to perform better for what concerns accuracy and a lower number of training snapshots.
\\
Possible future investigations could concern applications to problems with more involved stochastic dependence, as well as non-affinely parametrized problems. The latter ones could require the use of an ad-hoc weighted empirical interpolation technique \cite{peng:eim}. Another extension, especially in the greedy case, would be that of providing accurate estimation for the error. Such estimation were obtained for linear elliptic coercive problems in \cite{peng:w_elliptic_rnd}, but it would be useful to generalize them to different problems.
Finally, the proposed tests and methodology could also be used as the first step to study non-linear problems.

\paragraph{Acknowledgements}
We acknowledge the support by European Union Funding for Research and Innovation - Horizon 2020 Program - in the framework of European Research Council Executive Agency: H2020 ERC Consolidator Grant 2015 AROMA-CFD project 681447 ``Advanced Reduced Order Methods with Applications in Computational Fluid Dynamics''. We also acknowledge the INDAM-GNCS projects ``Metodi numerici avanzati combinati con tecniche di riduzione computazionale per PDEs parame- trizzate e applicazioni'' and ``Numerical methods for model order reduction of PDEs''.
The computations in this work have been performed with RBniCS \cite{rbnics:link} library, developed at SISSA mathLab, which is an implementation in FEniCS \cite{logg2012automated} of several reduced order modelling techniques; we acknowledge developers and contributors to both libraries.

\bibliographystyle{plain}      
\bibliography{references}   

\end{document}